\documentclass[journal]{IEEEtran}

\usepackage{amsfonts}
\usepackage{graphicx}
\usepackage{epsfig,rotating}
\usepackage{amssymb,amsmath}
\usepackage{latexsym}
\usepackage{array,multirow}
\usepackage{booktabs}

\usepackage[T1]{fontenc}
\usepackage{amsmath}

\newcommand{\be}{\begin{equation}}
\newcommand{\ee}{\end{equation}}
\newcommand{\beaa}{\begin{eqnarray*}}
\newcommand{\eeaa}{\end{eqnarray*}}
\newcommand{\bea}{\begin{eqnarray}}
\newcommand{\eea}{\end{eqnarray}}
\newcommand{\lbl}{\label}

\newcommand{\ml}{\mathcal}

\newtheorem{theorem}{Theorem}
\newtheorem{prop}{Proposition}
\newtheorem{lemma}{Lemma}

\begin{document}

\title {SURE Information Criteria for Large Covariance Matrix Estimation and Their Asymptotic Properties}
\author{Danning~Li and~Hui~Zou\thanks{D. Li is Assistant Professor in School of Mathematics at Jilin University, email:  danningli@jlu.edu.cn. This work was completed when Li was a research associate in the Statistical Laboratory, University of Cambridge, Cambridge, CB3 0WB.}  \thanks{H. Zou is Professor of Statistics at University of Minnesota, e-mail: zouxx019@umn.edu.
Zou's research is partially supported by NSF grant DMS 1505111.} \thanks{Manuscript received October 4 2014,; revised October 29 2015.}
}

\maketitle

\begin{abstract}
Consider $n$ independent and identically distributed $p$-dimensional Gaussian random vectors with covariance matrix $\Sigma.$ The problem of estimating $\Sigma$ when $p$ is much larger than $n$ has received a lot of attention in recent years. Yet little is known about the information criterion for covariance matrix estimation. How to properly define such a criterion and what are the statistical properties? We attempt to answer these questions in the present paper by focusing on the estimation of bandable covariance matrices when $p>n$ but $\log(p)=o(n)$. Motivated by the deep connection between Stein's unbiased risk estimation (SURE) and  Akaike information criterion (AIC) in regression models, we propose a family of generalized SURE ($\text{SURE}_c$) indexed by $c$ for covariance matrix estimation, where $c$ is some constant. When $c$ is 2, $\text{SURE}_2$ provides an unbiased estimator of the Frobenius risk of the covariance matrix estimator. Furthermore, we show that by minimizing $\text{SURE}_2$ over all possible banding covariance matrix estimators we attain the minimax optimal rate of convergence under Frobenius norm and the resulting estimator behaves like the covariance matrix estimator obtained by the so-called oracle tuning. When the true covariance matrix is exactly banded, we prove that by minimizing $\text{SURE}_{\log(n)}$ we select the true bandwidth with probability tending to one. Therefore, our analysis indicates that $\text{SURE}_2$ and $\text{SURE}_{\log(n)}$ can be regarded as the AIC and Bayesian information criterion (BIC) for large covariance matrix estimation, respectively.
\end{abstract}

\IEEEpeerreviewmaketitle

\section{Introduction}\lbl{intro}

Let $X_1,\cdots,X_n$ be independent and identically distributed $p$-dimensional Gaussian random vectors with mean $\mu$ and  covariance matrix $\Sigma_{p\times p}=\big(\sigma_{ij}\big)_{p\times p}$.  Recently the problem of large covariance matrix estimation has become a hot topic, because the classical sample covariance matrix performs very poorly in the high-dimensional setting (Johnstone, 2001). Several regularized estimators of large covariance matrices have been proposed, including banding (Wu and Pourahmadi, 2003; Bickel and Levina, 2008a; Xiao and Bunea 2014), tapering (Furrer and Bengtsson, 2007; Cai, Zhang and Zhou, 2010) and thresholding (Bickel and Levina, 2008b; El Karoui, 2008; Rothman, Levina and Zhu, 2009;  Cai and Liu 2011). In particular, several papers have been devoted to the study of bandable covariance matrices where the off-diagonal elements decay as they move away from the diagonal. Bandable covariance matrices have natural applications for modeling temporal and spatial dependence. Cai, Zhang and Zhou (2010) developed the first minimax optimality theory for estimating large bandable matrices. Consider the following  parameter space
\begin{eqnarray}\label{cai}
\mathcal{F}_{\alpha}&=&\{\Sigma: \max_{j}\sum_{i}\{|\sigma_{ij}|: |i-j|>k\}\leq Mk^{-\alpha} \,\\
&&\text{for} \, \text{all}\,  k, \text{and}\,  \lambda_{\max}(\Sigma)\leq M_0 \}, \nonumber
\end{eqnarray}
where $\lambda_{\max}(\Sigma)$ is the largest eigenvalue of  matrix $\Sigma$,  $\alpha, M, M_0$ are positive constants. The constant $\alpha $ controls the decay rate of the off-diagonal elements of $\Sigma$.
Cai, Zhang and Zhou (2010) established the minimax optimal rates of estimating $\Sigma$ under matrix $\ell_2$ and Frobenius norms and showed that tapering the sample covariance matrix with different tapering parameters can attain the minimax rates under matrix $\ell_2$ and Frobenius norms. Bien, Bunea and Xiao (2014) proposed a new convex banding estimator defined as the minimizer of a convex objective function and proved its minimax optimality. Qiu and Chen (2012) studied how to test whether a covariance matrix is banded or not.

In our work we consider a so-called generalized tapering estimator that includes the banding estimator in Bickel and Levina (2008a) and the tapering estimator in Cai, Zhang and Zhou (2010) as special cases. The performance of the generalized tapering estimator depends on the choice of the tapering parameter. From the application perspective, the most common practice is choosing a good tapering parameter based on the observed data. Given a target unknown covariance matrix and a series of tapering covariance estimators, the ideal selection method should choose the best tapering covariance estimator among all candidate tapering estimators. To handle the selection problem we need an objective criterion for comparing different covariance matrix estimators. It is well known that for comparing different regression models one can apply various information criteria such as AIC or BIC. The best selected regression model is the one with the smallest AIC or BIC score. The choice of the information criterion depends on the user's objective: AIC is used for optimizing the prediction performance while BIC is used to select the true model (Shao, 1997; Yang 2005). Given the new covariance matrix estimation problem and the great success of information criteria for model selection, we naturally ask the following question: what is the analogue of AIC or BIC for the covariance matrix estimation problem? In this paper we provide a good answer to the question by focusing on the bandable covariance matrix estimation problem. Our solution is based on Stein's unbiased risk estimation theory (Stein, 1981; Efron 1986, 2004). Consider a tapering covariance matrix estimator $\hat{\Sigma}^{(\tau)}$ where $\tau$ denotes the tapering parameter. Let $\tilde{\Sigma}^s =\big(\tilde{\sigma}^s_{ij}\big)_{p\times p}$ be the sample covariance matrix.  It can be shown that
the Frobenius risk can be written as
 \begin{equation*}\lbl{ch}R (\tau)= \mathbb{E} \|\hat{\Sigma}^{(\tau)}-\tilde{\Sigma}^s\|_{F}^2-\sum_{i,j}\text{var}(\tilde{\sigma}^s_{ij})+2\sum_{i,j} \text{cov}(\hat{\sigma}_{ij}^{(\tau)}, \tilde{\sigma}^s_{ij}).\end{equation*}
We can further derive $\text{SURE}(\tau)$, a Stein's unbiased risk estimator of $R(\tau)$, such that $\mathbb{E} [\text{SURE}(\tau)] = R(\tau)$.
Efron showed that for regression models with an additive homoscedastic Gaussian noise SURE is equivalent to AIC (Efron 1986, 2004). Motivated by this deep connection, we could regard the SURE formula as AIC for bandable covariance matrices. Furthermore, we extend SURE to a family of generalized SURE formulae denoted by $\text{SURE}_c$ where $c$ is a constant that may depend on $n,p$. Such an extension is motivated by the connections among different information criteria.
AIC is one of the information criteria defined as $$-2 \textrm{log likeligood}+c \times \textrm{model size}.$$ AIC uses $c=2$ while BIC uses $c=\log (n)$. We define $\text{SURE}_c(\tau)$ and $R_c (\tau)$ such that $\mathbb{E} [\text{SURE}_c(\tau)] = R_c(\tau)$ where
\begin{equation*}\lbl{chc}R_c (\tau)= \mathbb{E} \|\hat{\Sigma}^{(\tau)}-\tilde{\Sigma}^s\|_{F}^2-\sum_{i,j} \text{var}(\tilde{\sigma}^s_{ij})+c \sum_{i,j} \text{cov}(\hat{\sigma}_{ij}^{(\tau)}, \tilde{\sigma}^s_{ij}). \end{equation*}
The details of $\text{SURE}_c(\tau)$  are given in Section 2. We interpret $\text{SURE}_c$ as SURE information criteria for large covariance matrix selection: $\text{SURE}_2$ corresponds to AIC and $\text{SURE}_{\log(n)}$ corresponds to BIC.
In the sequel we reserve $\text{SURE}$ for $\text{SURE}_2$ to honor the literature on Stein's unbiased risk estimation.

We then study the asymptotic properties of $\text{SURE}_c$ and $\text{SURE}_c$ selection under the high-dimensional setting where $\log(p) =o(n)$, $p \ge n$ and $n \rightarrow \infty$. We make three theoretical contributions in this paper. First, we establish the central limit theorem of $\text{SURE}_c$ and a Bernstein type probability bound for $\text{SURE}_c$.  Next, we study the risk property of the SURE tuned estimator. We show that the SURE tuned estimator is minimax rate optimal under the Frobenius norm for estimating the parameter space $ \mathcal{F}_{\alpha}$ in  (\ref{cai}). This result holds for the generalized tapering estimator.  Moreover, we prove that the ratio of the risk of SURE tuned banding estimator to the smallest possible Frobenius risk by banding converges to 1 uniformly, which means that SURE tuning works as well as the oracle tuning. Last, we study the selection property of $\text{SURE}_c$ tuning when $c=2$ and $c=\log(n)$. To take advantage of its simplicity, we focus on the banding estimator in the study. Assume that $\Sigma$ is a banded matrix with bandwidth $k_0$. Under some regularity conditions, we show that the minimizer of SURE is in $[k_0, k_0+\log n]$ almost surely and the minimizer of $\text{SURE}_{\log(n)}$ equals $k_0$ almost surely. In summary, our theoretical results show that $\text{SURE}$ and $\text{SURE}_{\log(n)}$ have the fundamental properties of AIC and BIC (Shao, 1997; Yang 2005), that is, SURE is risk optimal and $\text{SURE}_{\log(n)}$ is selection consistent when the true model is in the candidate list (i.e., the true covariance matrix is banded).

The rest of the paper is organized as follows. Section 2 contains the discussion of SURE and its generalization to $\text{SURE}_c$ by following the information criteria argument. We prove the asymptotic normality of $\text{SURE}_c(\tau)-R_c(\tau)$ and further present a large deviation bound on  $| \text{SURE}_c(\tau)-R_c(\tau)|$. In Section 3 we study the risk property of the SURE tuned estimator and the selection property of SURE tuning and $\text{SURE}_{\log(n)}$ tuning. We conduct a simulation study to examine the theoretical results in Section 4. The proofs of main theorems are given in Section 5. We move the proofs of technical lemmas used in Section 5 to an appendix.

 \section{Limiting Results for SURE Information Criteria}
\subsection{$\text{SURE}$ information criteria}\lbl{Meili}

Let $X_1,\cdots,X_n$ be i.i.d. $p$-dimensional Gaussian random vectors with mean $\mu$ and  covariance matrix $\Sigma_{p\times p}=\big(\sigma_{ij}\big)_{p\times p}$. Let $\bar{X}=\frac{1}{n}\sum_{k=1}^nX_k$. The maximum likelihood estimation(MLE) of $\Sigma$ is
 $\tilde{\Sigma}=(\tilde{\sigma}_{ij})_{p\times p}=\frac{1}{n}\sum_{k=1}^n(X_k-\bar{X})(X_k-\bar{X})^T.$
We assume $p=p_n\geq n\geq 3$ and $\log p_n= o(n)$ in the sequel. It is assumed that $\Sigma$ comes from the parameter space $\mathcal{F}_{\alpha}$ defined  in $(\ref{cai})$.
Banding or tapering is a useful regularization method for estimating such covariance matrices (Bickel and Levina, 2008a; Cai, Zhang and Zhou, 2010).
To provide a unified treatment of banding and tapering, we consider the generalized tapering estimator of the covariance matrix:
\begin{equation*}
\hat{\Sigma}^{(\tau)}=\big(\hat{\sigma}_{ij}^{(\tau)}\big)_{1\leq i, j \leq p}=\big(\omega_{ij}^{(\tau)}\tilde{\sigma}_{ij}\big)_{1\leq i, j \leq p}
\end{equation*}
where   the generic tapering weights  $(\omega_{ij}^{(\tau)})_{1\leq i, j \leq p}
$ satisfy
  \begin{description}
 \item[(i)] $\omega_{ij}^{(\tau)}=1$ for $|i-j|\leq \lfloor \frac{\tau}{2}\rfloor$,
 \item[(ii)]$\omega_{ij}^{(\tau)}=0$ for $|i-j|\geq  \tau$,
  \item[(iii)] $0\leq \omega_{ij}^{(\tau)}\leq 1$ for $\lfloor \frac{\tau}{2}\rfloor<|i-j|< \tau$.
  \end{description}
It can be easily checked that both banding and tapering use some special weights and their weights satisfy  conditions (i)-(iii).  In some theorems we  consider the generalized tapering estimator, because the exact form of $ \omega_{ij}^{(\tau)}$ for  $\lfloor \frac{\tau}{2}\rfloor<|i-j|< \tau$ does not matter. In some theorems we focus on the banding estimator to take advantage of its simpler expression.

Let $R(\tau)=\mathbb{E} \|\hat{\Sigma}^{(\tau)}-\Sigma \|_{F}^2$ be the  Frobenius risk of $\hat{\Sigma}^{(\tau)}$. Yi and Zou (2013) proved the following Stein's identity
\begin{eqnarray*}
&&R(\tau)\\
&=&\mathbb{E} \|\hat{\Sigma}^{(\tau)}-\tilde{\Sigma}^s\|_{F}^2-\sum_{i,j}\text{var}(\tilde{\sigma}^s_{ij})+2\sum_{i,j} \text{cov}(\hat{\sigma}_{ij}^{(\tau)}, \tilde{\sigma}^s_{ij}) \nonumber  \\
&=&\mathbb{E} \|\hat{\Sigma}^{(\tau)}-\tilde{\Sigma}^s\|_{F}^2-\sum_{i,j}\text{var}(\tilde{\sigma}^s_{ij}) +2\frac{n-1}{n}\sum_{i,j} \omega_{ij}^{(\tau)}  \text{var} (\tilde{\sigma}^s_{ij}),
\end{eqnarray*}
where $\tilde{\Sigma}^s=\frac{n}{n-1}\tilde{\Sigma}$ is the sample covariance matrix of $\Sigma$ which is unbiased for $\Sigma$. The third term in the righthand is referred to as the covariance penalty (Efron, 2004).
Let $\widehat{\text{var}}(\tilde{\sigma}^s_{ij})$ be an unbiased estimator of $\text{var}(\tilde{\sigma}^s_{ij})$, then we have an unbiased estimator of $R(\tau)$ as follows (Yi and Zou, 2013)
\begin{eqnarray}\lbl{gen0}
\textrm{SURE}(\tau)&=&\|\hat{\Sigma}^{(\tau)}-\tilde{\Sigma}^s\|_{F}^2-\sum_{i,j} \widehat{\text{var}}(\tilde{\sigma}^s_{ij}) \\ \nonumber
&&+2 \frac{n-1}{n}\sum_{i,j} \omega_{ij}^{(\tau)} \widehat{\text{var}}(\tilde{\sigma}^s_{ij}).
\end{eqnarray}
Moreover, $\widehat{\text{var}}(\tilde{\sigma}^s_{ij})$ has an explicit expression and $\textrm{SURE}(\tau)$ is equal to
\begin{eqnarray}\lbl{gen1}
&&\textrm{SURE}(\tau)\\ \nonumber
&=&\sum_{1\leq i, j \leq p}(\frac{n}{n-1}-\omega_{ij} ^{(\tau)})^2\tilde{\sigma}_{ij}^2\\ \nonumber
&&+\sum_{1\leq i, j \leq p}(2\omega_{ij} ^{(\tau)}-\frac{n}{n-1})(a_n\tilde{\sigma}_{ij}^2+b_n\tilde{\sigma}_{ii}\tilde{\sigma}_{jj})\nonumber
\end{eqnarray}
with\begin{equation*} a_n=\frac{n(n-3)}{(n-1)(n-2)(n+1)}\, \text{ and}\, \, b_n=\frac{n}{(n+1)(n-2)}.\lbl{bed}\end{equation*}
$\textrm{SURE}(\tau)$ is referred to as Stein's unbiased risk estimator of $R(\tau)$. With SURE, one can select the tapering parameter by SURE tuning
\begin{equation}\lbl{coach0}
\hat \tau_n=\arg\min _{\tau}\text{SURE} (\tau).
\end{equation}
Prior to the covariance matrix estimation, SURE and SURE tuning have been used in SureShrink for adaptive wavelet thresholding (Donoho and Johnstone, 1995) and SURE-Lasso for high-dimensional linear model selection (Efron et al., 2004; Zou, Hastie and Tibshirani, 2007).

There is a deep connection between SURE and AIC in the context of regression analysis (Efron, 1986, 2004): the two are identical when the regression model has an additive homoscadestic Gaussian noise with known variance.
If fact, it has been argued that the covariance penalty in the SURE formula should be a universal way to define the degrees of freedom of an estimator (Efron, 1986, 2004; Zou, Hastie and Tibshirani, 2007). It is now well known that AIC is one of the many information model selection criteria defined as $-2 \textrm{log likeligood}+c \times \textrm{model size}$. AIC uses $c=2$ while BIC uses $c=\log (\textrm{sample size})$. It is interesting to see that the constant 2 appears naturally in the covariance penalty term in the SURE formula.  If we view $\textrm{SURE}(\tau)$ as the matrix counterpart of AIC, it is natural to ask what is the matrix counterpart of BIC? Motivated by the expression of information model selection criteria,  we replace the constant 2 in the covariance with $c$, while $2 \leq c=o(n)$.  Thus we define the generalized SURE formula as follows
\begin{eqnarray}\lbl{gen}
&&\text{SURE}_c(\tau)    \\
&=&\|\hat{\Sigma}^{(\tau)}-\tilde{\Sigma}^s\|_{F}^2-\sum_{i,j} \widehat{\text{var}}(\tilde{\sigma}^s_{ij})+c\frac{n-1}{n} \sum_{i,j} \omega_{ij}^{(\tau)} \widehat{\text{var}}(\tilde{\sigma}^s_{ij})\nonumber \\
&=&\sum_{1\leq i, j \leq p}(\frac{n}{n-1}-\omega_{ij} ^{(\tau)})^2\tilde{\sigma}_{ij}^2 \nonumber \\
&&+\sum_{1\leq i, j \leq p}(c \omega_{ij} ^{(\tau)}-\frac{n}{n-1})(a_n\tilde{\sigma}_{ij}^2+b_n\tilde{\sigma}_{ii}\tilde{\sigma}_{jj}).\nonumber
\end{eqnarray}
We define a new risk function $R_c(\tau)$ as $$R_c(\tau)=\mathbb{E} [\text{SURE}_c(\tau)],$$ then by (\ref{gen0}), (\ref{gen1}) and (\ref{gen}) we have
\begin{eqnarray}\lbl{BIC}
&&R_c(\tau) \\
&=& \mathbb{E} \|\hat{\Sigma}^{(\tau)}-\tilde{\Sigma}^s\|_{F}^2-\sum_{i,j}\text{var}(\tilde{\sigma}^s_{ij})+c\sum_{i, j} \text{cov}(\hat{\sigma}_{ij}^{(\tau)}, \tilde{\sigma}^s_{ij}). \nonumber
\end{eqnarray}
Naturally, we consider the minimizer of $\text{SURE}_c(\tau)$ as a chosen tapering parameter:
 \begin{equation*}\lbl{coach}
\hat{\tau}_{n}^{c}=\arg \min_{\tau}\text{SURE}_c(\tau).
  \end{equation*}
When $c=2$, (\ref{coach}) reduces to (\ref{coach0}). When $c=\log(n)$, we interpret (\ref{coach}) as BIC tuning. This interpretation will be rigorously  justified later in this paper.
Thus we treat $\text{SURE}_c(\tau)$ as a family of SURE information criteria for large covariance matrix estimation.

\subsection{The central limit theorem and large deviation bounds}
In this section, we establish the asymptotic distributional properties of the generalized SURE under the setting $p=p_n\ge n \rightarrow \infty$. Since we give a unified treatment of all $\text{SURE}_c$ we also write $c=C_n$ in the sequel to indicate the possible dependence of $c$ on $n$. For the sake of completeness, we restate the  following two assumptions:
 \begin{description}
 \item[(C.1)]  $X_1,\cdots, X_{n} $ are independent and identically distributed   $p$-dimensional Gaussian random vectors with mean $\mu$ and covariance $\Sigma$.
 \item[(C.2)] The covariance $\Sigma$ comes from the parameter space $\mathcal{F}_{\alpha}$ defined  in (\ref{cai}).
  \end{description}

To facilitate the analysis, we represent the $\text{SURE}_c(\tau)$ formula  by a new identity. According to Theorem 3.1.2 from Muirhead (1983), the MLE $\tilde{\Sigma} $ has the same distribution as $\hat{\Sigma}=\big(\hat{\sigma}_{ij}\big)_{1\leq i, j \leq p}=\frac{1}{n}\sum_{k=1}^{n-1}Z_kZ_k^{T}$, where $Z_j=(z_{ji})_{1\leq i\leq p},1 \le j \le n-1$ are i.i.d.  $N_{p}(0, \Sigma )$ random vectors.  Define \begin{eqnarray*}\lbl{devil}
S_n^c(\tau)=\sum_{1\leq i, j \leq p}\bar{a}_{ij}^{(\tau)}\hat{\sigma}_{ij}^2+\sum_{1\leq i, j \leq p}\bar{b}_{ij}^{(\tau)}(a_n\hat{\sigma}_{ij}^2+b_n\hat{\sigma}_{ii}\hat{\sigma}_{jj}),
\end{eqnarray*}
where
\begin{eqnarray*}
\bar{a}_{ij}^{(\tau)}=(\frac{n}{n-1}-\omega_{ij} ^{(\tau)})^2\,\text{and}\,\, \bar{b}_{ij}^{(\tau)}=(C_n\omega_{ij} ^{(\tau)}-\frac{n}{n-1}).\lbl{care}
 \end{eqnarray*}
Then we can conclude that $S_n^c(\tau)$ has the same distribution as $\text{SURE}_c(\tau)$. Therefore, it suffices to investigate the distributional properties of $S_n^c(\tau)$.

We derive a decomposition of $S_n^{c}(\tau)-R_c(\tau)$. Note that $R_c(\tau)=\mathbb{E} S_n^c(\tau)$. Define
$\mu_{ij}= \mathbb{E}\hat{\sigma}_{ij}^2=\frac{n(n-1)}{n^2}\sigma_{ij}^2+\frac{n-1}{n^2}\sigma_{ii}\sigma_{jj}$ and $\mu_{ij}'= \mathbb{E}\hat{\sigma}_{ii}\hat{\sigma}_{jj}=\frac{(n-1)^2}{n^2} \sigma_{ii}\sigma_{jj}+\frac{2(n-1)}{n^2}\sigma_{ij}^2$,
then we have
  \begin{eqnarray}\lbl{may}
 &&S_n^c(\tau)-R_c(\tau)\\
&=&\sum_{1\leq i, j \leq p}(\bar{a}^{(\tau)}_{ij}+a_n\bar{b}^{(\tau)}_{ij})(\hat{\sigma}_{ij}^2-\mu_{ij}) \nonumber \\
&&+\sum_{1\leq i, j \leq p}b_n\bar{b}^{(\tau)}_{ij}(\hat{\sigma}_{ii}\hat{\sigma}_{jj}-\mu_{ij}')\nonumber\\
 &=& \sum_{m=2}^{ n-1}\sum_{l=1}^{m-1}H_n(Z_m, Z_l)+\sum_{m=1}^{n-1}Y_m+\sum_{m=1}^{n-1}U_m +R_1+R_2,\nonumber
 \end{eqnarray}
 where
 \begin{equation}\lbl{nice}
 \begin{aligned}
  &H_n(Z_m, Z_l) \\
&=\frac{1}{n^2}\sum_{1\leq i, j \leq p}2\bar{A}^ {(\tau)}_ {ij}(z_{mi}z_{mj}-\sigma_{ij})(z_{li}z_{lj}-\sigma_{ij}),  \\
 Y_m&= \frac{2(n-2)}{n^2} \sum_{1\leq i, j \leq p}\bar{A}^ {(\tau)}_ {ij}\sigma_{ij} ( z_{mi} z_{mj}-\sigma_{ij}),\\
  U_m&= \frac{1}{n^2}  \sum_{1\leq i, j \leq p}\bar{B}^ {(\tau)}_ {ij}\{\sigma_{ii}(z_{mj}^2-\sigma_{jj})
 +\sigma_{jj}(z_{mi}^2-\sigma_{ii})\},\\
 R_1&= \frac{1}{n^2}\sum_{m=1}^{n-1}\sum_{1\leq i, j \leq p}
 \bar{C}^{(\tau)}_ {ij} \{( z_{mi}^2-\sigma_{ii})( z_{mj}^2-\sigma_{jj})-2\sigma_{ij}^2\},\\
     R_2&=  \frac{1}{n^2}\sum_{m=2}^{n-1}\sum_{l=1}^{ m-1}\sum_{1\leq i, j \leq p}
     2b_n\bar{b}^{(\tau)}_{ij}(z_{mi}^2-\sigma_{ii})(z_{lj}^2-\sigma_{jj})
     \end{aligned}
  \end{equation}
with
\begin{equation}\lbl{nice0}
 \begin{aligned}
\bar{A}^{(\tau)}_ {ij}&=\bar{a}_{ij}^{(\tau)}+a_n\bar{b}_{ij}^{(\tau)},\\
 \bar{B}^{(\tau)}_ {ij}&=\bar{a}_{ij}^{(\tau)}+(a_{n}+b_n (n-1))\bar{b}_{ij}^{(\tau)} ,\\
\bar{C}^{(\tau)}_ {ij}&=\bar{a}_{ij}^{(\tau)}+(a_{n}+b_n) \bar{b}_{ij}^{(\tau)}.
     \end{aligned}
  \end{equation}

We are interested in the asymptotic distribution of $\text{Var}_n(\tau)^{-\frac{1}{2}} (\text{SURE}_c(\tau)-R_c(\tau))$
 where $\text{Var}_n(\tau)$ is defined below
\begin{eqnarray}\lbl{sam} && \sum_{1\leq i, j ,s,t \leq p}\bigg\{ \frac{2(n-2)}{n^4}\bar{B}^ {(\tau)}_ {ij}\bar{B}^ {(\tau)}_ {st}(\sigma_{ii}\sigma_{ss}\sigma_{jt}^2+\sigma_{ii}\sigma_{tt}\sigma_{js}^2 \nonumber \\ &&+\sigma_{jj}\sigma_{ss}\sigma_{it}^2+\sigma_{jj}\sigma_{tt}\sigma_{is}^2) \nonumber \\
&& +
\frac{2(n-1)(n-2)}{n^4}\bar{A}^ {(\tau)}_ {ij}\bar{A}^ {(\tau)}_ {st} (\sigma_{is}\sigma_{jt}+\sigma_{it}\sigma_{js})^2 \nonumber \\
&&+ \frac{4(n-2)^3}{n^4}\bar{A}^ {(\tau)}_ {ij}\bar{A}^ {(\tau)}_ {st}\sigma_{ij}\sigma_{st} \times(\sigma_{is}\sigma_{jt}+\sigma_{it}\sigma_{js})\nonumber \\
&&+\frac{8(n-2)^2 }{n^4}\bar{A}^ {(\tau)}_ {ij}\bar{B}^ {(\tau)}_ {st}\sigma_{ij}( \sigma_{ss}\sigma_{it} \sigma_{jt}+\sigma_{tt}\sigma_{is}\sigma_{js})\bigg\}.\lbl{var}
\end{eqnarray}
Actually, $\text{Var}_n(\tau)$ approximates variance of $\text{SURE}(\tau)-R_c(\tau),$ by deleting  the higher order  terms come from $R_1$ and $R_2$.
From the definitions of $\text{Var}_n(\tau)$ and $R_c(\tau)$, we have the following proposition.
\begin{prop}\lbl{Meng} For any $1\leq \tau\leq p$ and $C_n=o(n)$, there exists a constant $C$ such that $\sqrt{\mathrm{Var}_n(\tau)}/R_{c}(\tau)\leq C\max(\frac{1}{\tau C_n}, \frac{1}{(np)^{1/2}})$.
 \end{prop}

 Assume $p=p_n\geq n$ and $C_n=o(n)$. In the following two theorems we prove the asymptotic normality of $\text{Var}_n(\tau)^{-\frac{1}{2}} (\text{SURE}_c(\tau)-R_c(\tau))$ under three asymptotic settings: (1). $p_n/n\to \infty$, (2). $C_n\to \infty$, (3). $C_n$ is constant and  $p_n/n\to y \in[1, \infty)$.

\begin{theorem}\lbl{girl} Suppose $(C.1)$ and $(C.2)$ hold. Assume $p=p_n\geq n$ and $2\leq C_n=o(n)$,  then for any $1\leq \tau\leq p$,   $\mathrm{Var}_n(\tau)^{-\frac{1}{2}} (\mathrm{SURE}_c(\tau)-R_c(\tau)) $ converges to the standard normal distribution if (i) $n\to\infty$ and $p_n/n\to \infty$ or (ii) $C_n\to \infty$ as $n\to\infty$.
\end{theorem}

 \begin{theorem}\lbl{babygirl} Suppose $(C.1)$ and $(C.2)$ hold. Assume $C_n$ is a constant and  $p_n/n\to y \in[1, \infty)$ as $n\to \infty$, then for any $\epsilon>0$ and $1\leq \tau\leq (1-\epsilon)p$, $\mathrm{Var}_n(\tau)^{-\frac{1}{2}} (\mathrm{SURE}_c(\tau)-R_c(\tau))$  asymptotically follows the standard normal distribution.\end{theorem}

In the next theorem we derive an explicit probability bound to describe how $\text{SURE}_c(\tau )$ deviates from $R_c(\tau)$.
         \begin{theorem}\lbl{zhen}  Suppose $(C.1)$ and $(C.2)$ hold. Assume $p=p_n\geq n$ and $2\leq C_n=o(n)$,
for any even number $K_0\geq 4$,      then there  exists $M_{K_0}$ such that, for any $\epsilon>0$ and $1\leq \tau\leq(1-\epsilon)p$,   $ \mathbb{P}(  |\mathrm{SURE}_c(\tau)-R_c(\tau) |\geq \lambda t_n \sqrt{\mathrm{Var}_n(\tau)} )\leq  2 \exp\big(-\frac{1}{16}t_n^2\big)+    M_{K_0}\big\{(np)^{-\frac{K_0}{2}} + + (\lambda^2 t_n^2n)^{-K_0/2}+  \frac{ t_n^{K_0-4 }}{\lambda ^{K_0 }n^{K_0/2-1}}+\frac{C_n^{K_0}}{(\lambda t_n)^{K_0}n^{ K_0}}\big \}$ for all  $\lambda\geq 1$, $t_n>0$ and $n\geq 3$. \end{theorem}

Combining Proposition \ref{Meng} with Theorem~\ref{zhen}, we easily have the following proposition.

      \begin{prop}  $\mathrm{SURE}_c(\tau)/R_c( \tau) \to 1$ in probability when $C_n\to \infty$ as $n \to \infty$. If $2\leq C_n$ is a constant, then $\mathrm{SURE}(\tau_n)/R( \tau_n)\to 1$ in probability as $\tau_n\to \infty$ when $n\to \infty$. \end{prop}

    \section{Properties of The $\text{SURE}_c$ Tuned Estimators}
In this section we study the asymptotic properties of the  $\text{SURE}_c$ tuned estimator and  $\text{SURE}_c$ selection.
To honor the literature we reserve SURE for $\text{SURE}_2$.

\subsection{Minimax optimality of SURE}

AIC is known to yield an asymptotic minimax estimator (Yang, 2005). We interpret SURE as the AIC for covariance matrix estimation. Thus we expect the same minimax optimality property holds
for SURE.

\begin{theorem}\lbl{cooked} Suppose $(C.1)$ and $(C.2)$ hold.
Assume $n\leq p $ and $\log p =o(n)$, then $\sup\limits_{\Sigma\in {\cal F}_{\alpha}}\mathbb{E}\|\hat{\Sigma}^{(\hat{\tau}_n)}-\Sigma\|_F^2\asymp pn^{-(2\alpha+1)/2(\alpha+1)}$.
\end{theorem}

Cai, Zhang and Zhou (2010) showed that the minimax rate of convergence for estimating $\Sigma$ in ${\cal F}_{\alpha}$ under Frobenius norm is $pn^{-(2\alpha+1)/2(\alpha+1)}$.
Thus Theorem~\ref{cooked} indicates that SURE tuning yields a minimax rate optimal tapering estimator for estimating $\Sigma$ in ${\cal F}_{\alpha}$. The estimator defined in Cai, Zhang and Zhou (2010) that attains the minimax rate depends on knowing
$\alpha$. The SURE tuned estimator is fully data-driven. Cai and Yuan (2012) constructed another fully data-driven minimax rate optimal estimator by using the idea of block-thresholding.

\subsection{SURE tuning versus oracle tuning}
Define $\tau_0=\arg\min_{\tau} R(\tau)$. Then $\hat \Sigma^{(\tau_0)}$ is called the oracle tuned estimator because it yields the smallest risk. The oracle tuning only exists in theory but can be used to judge the performance of an actual tuning method.
We compare SURE tuning with the oracle tuning.

We begin with some regularity conditions.
We switch the parameter space from ${\cal F}_{\alpha}$ to a slightly different parameter space ${\cal G}_{\alpha}$ by following Cai, Zhang and Zhou (2010). Correspondingly, we replace the condition (C.2) by (C.3)
\begin{description}
 \item[(C.3)] $\Sigma$ is in $\mathcal{G}_{\alpha}$ where
\begin{eqnarray}
\mathcal{G}_{\alpha}&=&\{\Sigma:   |\sigma_{ij}|\leq M_1|i-j|^{-(\alpha+1)}\, \nonumber \\
&&\text{for} \, \text{all}\, i\neq j\, \text{and}\, \lambda_{max}(\Sigma)\leq M_0 \}.
\end{eqnarray}
   \end{description}
It is worth mentioning that estimating ${\cal G}_{\alpha}$ is as hard as estimating  ${\cal F}_{\alpha}$ because the minimax rate stays the same (Cai, Zhang and Zhou, 2010). We work with ${\cal G}_{\alpha}$ because it makes our analysis
slightly easier. For the same reason of convenience, we focus on the banding estimator instead of the generalized tapering estimator.
 \begin{description}
 \item[(C.4)]  The tapering weights $w_{ij}^{(\tau)}$ is $w_{ij}^{(\tau)}= I (|i-j|< \tau)$ for all $1\leq i, j\leq p $.
  \end{description}
 Under condition (C.4), we can have a simpler expression of $R(\tau)$:
\begin{eqnarray}\lbl{vapor} R(\tau) = \sum_{|i-j|< \tau}\big(\frac{1}{n}\sigma_{ij}^2+\frac{n-1}{n^2}\sigma_{ii}\sigma_{jj}\big)+\sum_{| i- j|\geq \tau}  \sigma_{ij}^2.\end{eqnarray}
 We make additional assumption on the covariance matrix $\Sigma$:
  \begin{description}
    \item[(C.5)] There exists $\gamma> 1$ such that $\sum_{|i-j|=k}\sigma^2_{ij}\geq \gamma \sum_{|i-j|=k+1}\sigma^2_{ij}$ for all $k\geq 0$ and $\sigma_{ii}=1$ for all $1\leq i\leq p$.
   \end{description}
 Condition (C.5) is not very strict. It only requires  the decay trend is detectable when  the covariances $\sigma_{ij}$ move away from the diagonal. The assumption of  $\sigma_{ii}=1$ is just used for simplifying the proof.

\begin{theorem}\lbl{daughter}
Suppose $(C.1)$ and $(C.3)-(C.5)$ hold. Let  $\tau_0 $ is the unique minimizer of $R(\tau)$ and $\hat{\tau}_{n}=\arg \min_{\tau}\mathrm{SURE}(\tau)$. Assume $n\leq p $ and $\log p =o (n)$, then $ |\hat{\tau}_n-\tau_0|\leq \log n $ almost surely as $n\to \infty$. Further assume that $n \log n \ll p$, then $\sup\limits_{\Sigma\in \mathcal{G}_{\alpha}}|\mathbb{E}\|\hat{\Sigma}^{(\hat{\tau}_n)}-\Sigma\|_{F}^2/R(\tau_0)- 1|\to 0$ as $n\to \infty.$
\end{theorem}

Theorem~\ref{daughter} shows that if we only care about the risk property of the estimator, SURE tuning works as well as the oracle tuning because the SURE tuned banding estimator automatically achieves the smallest Frobenius risk among all possible banding estimators. The same conclusion can be established for the tapering estimator proposed in Cai, Zhang and Zhou (2010).  Its proof is slightly more involved. For the sake of space we do not include it here.

\subsection{Bandwidth selection by SURE and $\text{SURE}_{\log(n)}$}

Theorem~\ref{daughter} also shows that the distance between the SURE selection result and $\tau_0$ is bounded by $\log(n)$.
When $\Sigma$ is a banded matrix, then we expect that $\tau_0$ is the bandwidth of $\Sigma$.
Thus, it is interesting to see if SURE selection could always pick $\tau_0$. Note that when using tapering weight, $R(\tau)$ reaches its minimum at $\tau_0= 2k_0-3$, where $k_0$ is the bandwidth of $\Sigma$. For simplicity, we consider the banding covariance matrix estimator under the assumption that the true covariance matrix is exactly banded.

\begin{theorem}\lbl{Dad}  Suppose $(C.1)$  and $(C.4)$ hold. Let the covariance matrix $
  \Sigma$ be a banded matrix with bandwidth $k_0$ such that $\sigma_{ij}=0$ if $|i-j|\geq k_0$ and   $\min\limits_{|i-j|\leq k_0-1}\sigma_{ij}^2\gg \log n/n$, where $k_0$ is a constant does not depend on $n$. If $n \leq p$ and  $\log p =o(n)$, then $k_0+\log n\geq  \hat{\tau}_n\geq k_0 $ almost surely.\end{theorem}

Recall that in the context of linear regression,  BIC is known for its selection consistency property if the true model is in the list of candidate models.
We define $\text{SURE}_{\log(n)}$ from SURE by following the relation between BIC and AIC.
Thus we expect $\text{SURE}_{\log(n)}$ tuning is selection consistent. We make this claim rigorous in the sequel.

Under condition (C.4), we can write
\begin{eqnarray}\lbl{vapor1}
&&R_c(\tau) \nonumber \\
&=&\sum_{|i-j|< \tau}\big\{\big(\frac{c}{n}-\frac{1}{n}\big)\sigma_{ij}^2+\big(\frac{n-1}{n^2}+\frac{c-2}{n}\big)\sigma_{ii}\sigma_{jj}\big\} \nonumber \\&&+\sum_{| i- j|\geq \tau}  \sigma_{ij}^2.
\end{eqnarray}
Define $\tau_0^c=\arg\min_{\tau} R_c(\tau)$.

\begin{theorem}\lbl{bake}
Suppose $(C.1)$ and $(C.3)-(C.5)$ hold. Assume that  $\tau_0^c$ is the unique minimizer of $R_c(\tau)$. Let $c= \log n$, if  $n\leq p$ and $\log p =o(n)$, then $ |\hat{\tau}_n^{c}-\tau_0^c|\leq 1 $ almost surely. If exists $\delta>0$ such that $|R_c(\tau_0^c\pm 1)-R_c(\tau_0^c)|\geq 2\delta \log n\sqrt{\mathrm{Var}_n(\tau_0^c)}$, then $ \hat{\tau}_n^{c}=\tau_0^c $ almost surely as $n\to \infty$.
\end{theorem}

Now we assume that $\Sigma$ is a banded matrix with bandwidth $k_0$. The next theorem shows that $\text{SURE}_{\log(n)}$ tuning selects the true bandwidth $k_0$ almost surely.
\begin{theorem}\lbl{breaking} Suppose $(C.1)$  and $(C.4)$ hold. Let the covariance matrix $
  \Sigma$ be a banded matrix with bandwidth $k_0$ such that $\sigma_{ij}=0$ if $|i-j|\geq k_0$ and   $\min\limits_{|i-j|\leq k_0-1}\sigma_{ij}^2\gg \log n/n$, where $k_0$ is a constant doesn't depend on $n$. If $n \leq p_n$ and  $\log p_n=o(n)$, let $c=\log n$ then $\hat{ \tau}_n^{c}=k_0 $ almost surely.\end{theorem}

\section{A Simulation Study}

In this section, we conduct a small simulation study to show $\text{SURE}$ selection is risk optimal and $\text{SURE}_{\log(n)}$ selection is consistent.

The simulated data were generated form $N(0,\Sigma)$ where three covariance models were considered.
\begin{itemize}
\item \texttt{Model 1} The covariance matrix has the form
\begin{equation*}\lbl{model1}
 \sigma_{ij}= \left\{
  \begin{array}{l l l}
 1,& 1\leq i=j\leq p\\   \rho|i-j|^{-(\alpha+1)}& 1\leq i\neq j\leq p.
 \end{array} \right.
\end{equation*}
We let $\rho=0.6$, $\alpha=0.1, 0.5,$ $n=250$ and $p=500$.

\item \texttt{Model 2} The covariance matrix has the form $\sigma_{ij}=\rho^{|i-j|}$, $1\leq i, j\leq p.$ We let $\rho=0.95, 0.5$, $n=250$ and $p=500$.

\item \texttt{Model 3}   The covariance matrix has the form $\sigma_{ij}=I(i=j)+\frac{1}{4}I(|i-j|\leq 4)$  $1\leq i, j\leq p.$  We let   $n=250$ and $p=500, 1000$.  This covariance is a banded  matrix with bandwidth $5$.
  \end{itemize}

We used the banding estimator to estimate $\Sigma$. Models 1 and 2 were used to test SURE selection and model 3 was used to test $\text{SURE}_{\log(n)}$ selection. Simulation results are summarized in Tables 1 and 2 where we report the average values based on 100 independent replications and the corresponding standard errors are shown in parenthesis.

Table 1 shows that SURE selection leads to the risk optimal estimator and Table 2 shows that $\text{SURE}_{\log(n)}$ selection is consistent in identifying the true bandwidth.

\begin{table}[ht]
\caption{Examine the optimal risk (Frobenius risk) property of SURE selection. The second column is the minimum of the Frobenius risk over all possible banding covariance matrix estimators. The third column is the Frobenius risk of the SURE tuned banding covariance matrix estimator.}
\centering
\begin{tabular}{c c  c }
\hline
Model 1& Minimum risk by Banding& \text{SURE} selected Banding\\
\hline
 {$\alpha=0.5$}& 58.83&57.57 (0.89)\\

{$\alpha=0.1$}& 30.16 &30.20 (0.67)\\
  \hline
Model 2& Minimum risk by Banding&\text{SURE} selected Banding\\
\hline
 {$\rho=0.95$}& 275.06 &273.48 (6.51)\\

 {$\rho=0.5$}&22.37 & 22.675 (0.71)\\
    \hline
\end{tabular}
\label{table:nonlin}
\end{table}

\begin{table}[ht]
\caption{Examine the selection consistency of $\text{SURE}_{\log(n)}$ selection. The second column is the true bandwidth of the covariance matrix in Model 3.
The third column is selected bandwidth by $\text{SURE}_{\log(n)}$ tuning.}
\centering
\centering
\begin{tabular}{c c  c }
\hline
Model 3 &true bandwidth&selected bandwidth \\
\hline
$p=500$  &  5 & 5 (0)\\
$p=1000$ & 5 &5 (0)\\
\hline
\end{tabular}
\label{table:nonlin}
\end{table}

\section{Proofs of Main Theorems}\lbl{proofs}In this section, we first list a few technical lemmas and then present the proofs of the main results.
The proofs of some technical lemmas are given in an appendix. Throughout this section we use $C$ to denote a generic constant.

 \subsection{Technical Lemmas}

The  general form of Isserlis' theorem, due to Withers (1985), is stated as follows.
\begin{lemma} \lbl{IS}If $A = \{ i_1,\cdots, i_{2N} \}$ is a set of integers such that $1
\leq i_k \leq p$, for any $k \in \{1,\cdots, 2N\}$ and $X=(X_i)_{1\leq i\leq p}\in  R^p$ is a Gaussian vector with zero mean then\begin{eqnarray*}
\mathbb{E}\prod_{i_k\in A}X_{i_k} =\sum\prod_{A }E(X_{i_l}X_{i_m}), \end{eqnarray*}
where $\sum\prod\limits_{A }$ denotes the sum of these products over all distinct ways pairing $\{1,\cdots, 2N\}.$
Moreover, if $A = \{ i_1,\cdots, i_{2N} , i_{2N+1}\}$ then, under the same assumptions, $\mathbb{E}\prod_{i_k\in A}X_{i_k}  = 0$.\end{lemma}

The following Bernstein-type inequality for martingale is essentially a special case of Theorem 1.2A  from   De La Pe$\tilde{n}$a (1999).
\begin{lemma}\lbl{DE} Let $\{d_i, \ml{F}_i\}$ is a martingale difference sequence with $\mathbb{E}(d_i|\ml{F}_{i-1})=0$ and $\mathbb{E}(d_i^2|\ml{F}_{i-1})=\sigma_{i}^2$, $V_n^2=\sum_{i=1}^n\sigma_{i}^2$. Furthermore, assume that $|d_i|\leq C$, for $0<C<\infty$. Then for any $x$, $y>0$
\begin{eqnarray*}
\mathbb{P}(\sum_{i}d_i\geq x, V_n^2\leq y)\leq \exp\bigg(-\frac{x^2}{2(y+Cx)}\bigg). \end{eqnarray*}
\end{lemma}

The following moment inequality for martingale  comes  from  Dharmadhikari et al. (1968).
\begin{lemma}\lbl{DATL} Let $\{d_i, \ml{F}_i\}$ is a martingale difference sequence with $\mathbb{E}(d_i|\ml{F}_{i-1})=0$, then for any $K_0\geq 2$ and $n\geq 1$, \begin{eqnarray*}
\mathbb{E}{|\sum_{i=1}^n d_i|^{K_0}}\leq L_{K_0}n^{K_0/2-1}\sum_{i=1}^{n}\mathbb{E}|d_i|^{K_0},
 \end{eqnarray*}where $L_{K_0}=\{8(K_0-1)\max(1, 2^{K_0-3})\}^{K_0}.$
\end{lemma}

From the paper by Cai et al. (2010), we get the following Lemma.
\begin{lemma}\lbl{Caizz}  Assume $\check{\Sigma}$ is an arbitrary estimator of covariance $\Sigma$. The minimax risk under the Frobenius norm satisfies that
\begin{eqnarray*}
&&\inf\limits_{\check{\Sigma}}\sup\limits_{\Sigma \in \mathcal{F}_{\alpha}}\mathbb{E}\|\check{\Sigma}-\Sigma\|_{F}^2\\
&\asymp&\inf\limits_{\check{\Sigma}}\sup\limits_{\Sigma \in\mathcal{G}_{\alpha}}\mathbb{E}\|\check{\Sigma}-\Sigma\|_{F}^2 \\
&\asymp &pn^{-(2\alpha+1)/2(\alpha+1)},
 \end{eqnarray*}
 Where $\mathcal{F}_{\alpha}$ and $\mathcal{G}_{\alpha}$ are defined as in (C.2) and (C.3).
  Furthermore for the parameter space $\mathcal{G}_{\alpha}$, the optimal tapering parameter $\tau$ of $\hat{\Sigma}^{(\tau)}$ is the order $n^{1/(2(\alpha+1))}$.
\end{lemma}

 Define $\mathcal{P}:=\{\{a_1, b_1\}, \{a_2, b_2\}, \cdots, \{a_{K_0}, b_{K_0}\} |\, a_i, b_i \in\{1,\cdots, 2K_0\} \}$ is a way of partitioning $\{1,\cdots, 2K_0\}$ into pairs.
 \begin{lemma}\lbl{iran}
         Suppose $Z=(z_{j})_{1\leq j\leq p}$  be a $p$-dimensional normal random vector with mean $0$ and covariance matrix $\Sigma_{p\times p}=\big(\sigma_{ij}\big)_{p\times p}$. Then $\mathbb{E}\prod_{k=1}^{K_0}(z_{i_{2k-1}}z_{i_{2k}}-\sigma_{i_{2k-1}i_{2k}})=\sum_{\mathcal{P}\in \mathcal{P}_{2K_0}}\prod_{\{s,t \}\in \mathcal{P}}\sigma_{i_si_t} $, where $K_0\geq 1$ and $\mathcal{P}_{2K_0}$ is  the set contains    all the distinct ways  that partitioning $\{1, \cdots,  2K_0\} $ into pairs excluding the ways partitioning $ 2k-1$ and $2k $ into a pair for any $k\leq K_0$. If $i_{2k-1}=i_{2k}$, for any even number $K_0\geq 2$ there exists C such that  $\sum_{1\leq i_1, \cdots, i_{K_0}\leq p}\mathbb{E}\prod_{k=1}^{K_0}(z_{i_{k}}^2-\sigma_{i_ki_{k}})\leq Cp^{K_0/2}$.
               \end{lemma}

       \begin{lemma}\lbl{rain}
Define $\Sigma_a=(|\sigma_{ij}|)$, where $\sigma_{ij}$
are the elements of $\Sigma$ and $\Sigma$ satisfies  the condition (C.2). Let $(i_1,\cdots, i_k)\in \{1, \cdots, p\}^k$, then
\begin{eqnarray*}
&&\sum_{1\leq i_1,\cdots, i_k\leq p}|\sigma_{i_1i_2}\sigma_{i_2i_3}\cdots \sigma_{i_ki_1}| \nonumber \\
&=&\mathrm{Trace}(\Sigma_a^k)\nonumber \\
&\leq &p(3M_0+M) ^k. \lbl{tian}
\end{eqnarray*}
\end{lemma}

   \begin{lemma}\lbl{day}
Consider  $\bar{A}^{(\tau)}_ {ij}$,   $\bar{B}^{(\tau)}_ {ij}$ and  $\bar{C}^{(\tau)}_ {ij}$ given in (\ref{nice0}). Assume $2\leq C_n=o(n)$ as $n\to \infty$, then $\bar{A}^{(\tau)}_ {ij}$and $\bar{C}^{(\tau)}_ {ij}$ are uniformly bounded by 2 for all $1\leq \tau\leq p$ when $n\geq 3$ . Furthermore $\bar{B}^{(\tau)}_ {ij}\leq C_n$ for all $\tau$, $\bar{B}^{(\tau)}_ {ij}=0$ for any $|i-j|\geq \tau$ and $\bar{B}^{(\tau)}_ {ij}
\geq \frac{C_n}{2}$ for any $|i-j|\leq \lfloor\frac{\tau}{2}\rfloor$.
\end{lemma}

Consider $H_n(Z_m, Z_l)$, $Y_m$, $U_m$, $R_1$ and $R_2$ defined in (\ref{nice}). We further define that $Y_{nm}=\sum_{l=1}^{m-1}H_n(Z_m, Z_l)$ for any $m\geq 2.$ and $Y_{nm}=0$ for $m=1$ and $G_m=Y_{nm}+Y_m+U_m.$  Then we have the following properties.
       \begin{lemma}\lbl{zhu}  Assume $p=p(n) \geq n$ and $C_n=o(n)$ as $n\to \infty$,
let
    \begin{eqnarray*} D_n&=&\sum_{m=2}^{n-1}\sum_{k=2}^{n-1} \mathbb{E}\{\mathbb{E}_{m-1}(Y_{nm}Y_m)\mathbb{E}_{k-1}(Y_{nk}Y_k) \\
    &&+\mathbb{E}_{m-1}(Y_{nm}U_m)\mathbb{E}_{k-1}(Y_{nk}U_k) \\
    &&+ \mathbb{E}_{m-1}(Y_{nm}Y_m)\mathbb{E}_{k-1}(Y_{nk}U_k)]\\
    &&+ \mathbb{E}_{m-1}Y_{nm}^2\mathbb{E}_{k-1}(Y_{nk}Y_k)\\
     &&+\mathbb{E}_{m-1}Y_{nm}^2\mathbb{E}_{k-1}(Y_{nk}U_k)\}.\end{eqnarray*}   Then there exists a constant C such that, for any $1\leq \tau\leq p$, $D_n\leq C(\frac{p(\tau C_n)^2}{n^5}+\frac{p^2 \tau C_n}{n^5}+\frac{p^2}{n^4})$.
  \end{lemma}

 \begin{lemma}\lbl{son}
     Let $Z_1,\cdots, Z_{n-1}$ be i.i.d. $p$-dimensional normal random vectors with mean $0$ and covariance matrix $\Sigma_{p\times p}=\big(\sigma_{ij}\big)_{p\times p}$. Assume $\Sigma$ satisfies (C.2). Let $s_{n-1}^2$ equal to $\mathbb{E}(\sum_{m=1}^{ n-1}G_m)^2$. Assume $p\geq n$ and $2\leq C_n=o(n)$ as $n\to \infty$, then $s_{n-1}^2$ has the form as in (\ref{var}) and there exists C such that, for any $1\leq \tau$ the upper bound $C(\frac{p^2}{n^2}+\frac{C_n^2}{n^3}  p\tau^2+\frac{C_n}{n^2} p\tau)$ and the lower bound is $\frac{1}{C} \max (\frac{p^2}{n^2},\frac{p\tau^2C_n^2}{n^3} )$ if (i) $p_n\gg n$ or (ii) $C_n$ depends on $n$ and $C_n\to \infty$ as $n\to\infty$.  \end{lemma}

    \begin{lemma}\lbl{night}  Let $Z_1,\cdots, Z_{n-1}$ be i.i.d. $p$-dimensional normal random vectors with mean $0$ and covariance matrix $\Sigma_{p\times p}=\big(\sigma_{ij}\big)_{p\times p}$. Assume $\Sigma$ satisfies  (C.2).
  Then for any even number $K_0$,   there exists a constant $C$ such that,  for any $1\leq \tau\leq p$, $ \mathbb{E}H_n(Z_1, Z_2)^{K_0} \leq  C\frac{p^{K_0}}{n^{2K_0}}$,   $\mathbb{E}Y_{m}^{K_0}\leq C\frac{p^{K_0/2}}{n^{K_0}}$ and    $\mathbb{E}U_m^{K_0}\leq C \frac{(C_n\tau)^{K_0}p^{K_0/2}}{n^{2K_0}}$  for all $n\geq 3$. \end{lemma}

  \begin{lemma}\lbl{dress}
 Let $Z_1,\cdots, Z_{n-1}$ be i.i.d. $p$-dimensional normal random vectors with mean $0$ and covariance matrix $\Sigma_{p\times p}=\big(\sigma_{ij}\big)_{p\times p}$. Assume $\Sigma$ satisfies  (C.2).    Then for any even number $K_0$, there exists a constant $C$ such that, for any $1\leq \tau\leq p$ and $n\geq 3$, $\mathbb{E}[\sum_{m=2}^{n-1}(\mathbb{E}_{m-1}Y_{nm}^2 -\mathbb{E}Y_{nm}^2)]^{K_0}\leq C(\frac{p^{3K_0/2}}{n^{5K_0/2}}+\frac{p^{K_0}}{n^{2K_0} })$.
\end{lemma}

 \begin{lemma}\lbl{night1} Let $Z_1,\cdots, Z_{n-1}$ be i.i.d. $p$-dimensional normal random vectors with mean $0$ and covariance matrix $\Sigma_{p\times p}=\big(\sigma_{ij}\big)_{p\times p}$. Assume $\Sigma$ satisfies  (C.2).
   Then  for any even number $K_0$, there exist a constant $C$  such that,  $1\leq \tau\leq p$ and $n\geq 3$,  $\mathbb{E}(R_1^{K_0})
\leq C\frac{p^{K_0}}{n^{3K_0/2}}$ and $ \mathbb{E}R_2^{K_0}\leq C\frac{(C_np)^{K_0}}{n^{2K_0}}$.
\end{lemma}

\begin{lemma}\lbl{fly}  Suppose $(C.1)$ holds and the covariance matrix $
  \Sigma$ is a banded matrix with bandwidth $k_0$ such that $\sigma_{ij}=a_{ij}I(|i-j|<k_0)$.   Further assume its nonzero elements are bounded away from zero, $|\sigma_{ij}|>b>0$ if $\sigma_{ij}\neq 0$ and $k_0=o(\sqrt{\frac{n}{C_n}})$. Consider $R_c(\tau)$ defined in (\ref{BIC}), where $2\leq c=C_n=o(n)$. With the tapering weight, $R_c(\tau)$ reaches its minimum at $2k_0-3$ for $k_0\geq 3$. With the banding weight, then $R_c(\tau)$ reaches its minimum at $k_0$.
\end{lemma}

\subsection{Proofs of Main Results}
\begin{IEEEproof}[Proof of Proposition \ref{Meng}]
By the definition of $R_{c}(\tau)$ and $\sigma_{ii}\geq c$ for all $1\leq i\leq p$, we know that there exists a constant $C_0$ such that
    \begin{eqnarray*}&&R_c(\tau)\\
    &=&\sum_{\lfloor \frac{\tau}{2}\rfloor<| i-j|< \tau} \{[\frac{n-1}{n}((\omega_{ij}^{(\tau)})^2-\frac{2n-C_n}{n-1}\omega_{ij}^{(\tau)})+1]\sigma_{ij}^2\\
    &&+\frac{n-1}{n^2}[(\omega_{ij}^{(\tau)})^2+\frac{n(C_n-2)}{n-1}\omega_{ij}^{(\tau)}]\sigma_{ii}\sigma_{jj}\}\\
    &&+\sum_{|i-j|\leq \lfloor \frac{\tau}{2}\rfloor}[(\frac{C_n-1}{n}\sigma_{ij}^2+\frac{nC_n-n-1}{n^2}  \sigma_{ii}\sigma_{jj}]\\
    &&+\sum_{| i- j|\geq \tau}  \sigma_{ij}^2\\
   &\geq& \frac{1}{C_0}\frac{p\tau C_n}{n}.
   \end{eqnarray*}
   Combining this with the upper bound of $\text{Var}_{n}(\tau)$ from Lemma \ref{son}, we conclude that
   $$\sqrt{\text{Var}_n(\tau)}/R_{c}(\tau)\leq C\max(\frac{1}{\tau C_n}, \frac{1}{(np)^{1/2}}).$$
\end{IEEEproof}

\begin{IEEEproof}[Proof of Theorem \ref{girl}] Since $S_n^c(\tau)$ has the same distribution as $\text{SURE}_c(\tau)$, we  consider $S_n^c(\tau)-R_c(\tau)$ instead of $\text{SURE}_c(\tau)-R_c(\tau)$.
Define
\begin{eqnarray*}\lbl{super}G_m=\sum_{l=1}^{m-1}H_n(Z_m, Z_l)+Y_m+U_m
\end{eqnarray*} for $m\geq 2$ and $G_1=Y_1+U_1$. Let  $S_{k}= \sum_{m=1}^{ k}G_m$ and
 $\ml{F}_k$ denote  the $\sigma$-field generated by $(Z_1,\cdots, Z_k)$ for $ k\geq 1$. 
 Then $\{S_n, \ml{F}_n, n=1,2,\cdots\}$ is a martingale on the probability space.  Notice that $S_n^c(\tau)-R_c(\tau)=S_{n-1}+R_1+R_2.$  First we  apply Brown's Martingale central limit theorem to $S_{n-1}$ (see Hall (1989) as an example).
Let  $\mathbb{E}_{k-1}$ denote $\mathbb{E}(\cdot |\ml{F}_{k-1})$,
\begin{eqnarray*} V_{n-1}^2=\sum_{m=1}^{n-1}\mathbb{E}_{m-1}G_m^2\, \text{and}\,   s^2_{n-1}=E(S_{n-1}^2)=E(V_{n-1}^2)\lbl{pain}.\end{eqnarray*}  We only need  to check the following  two conditions:
    \begin{eqnarray}\lbl{cat}
s_{n-1}^{-2}\sum_{m=1}^{n-1}   \mathbb{E}(G_m)^2I(|G_m|>\epsilon s_{n-1})\to 0,
    \end{eqnarray}
    as $n\to \infty$ for each $\epsilon>0,$ and
       \begin{eqnarray}\lbl{dog1}
       s_{n-1}^{-2}V_{n-1}^2\to 1\,\,\,\,\text{ in probability},
        \end{eqnarray}
         as $n\to \infty$. If the conditions  (\ref{cat}) and (\ref{dog1}) are satisfied, then $s_{n-1}^{-1}S_{n-1}$ is asymptotically $N(0,1).$

         First we prove (\ref{cat}).
           Define $Y_{nm}=\sum_{l=1}^{m-1}H_n(Z_m, Z_l)$ for $m\geq 2$ ,  using the $C_r$ inequality, we have
  \begin{eqnarray*}
\sum_{m=1}^{n-1} \mathbb{E} G_m^4 \leq 27\sum_{m=2}^{n-1} \mathbb{E}Y_{nm}^4 +27\sum_{m=1}^{n-1} \mathbb{E}Y_m^4+27\sum_{m=1}^{n-1} \mathbb{E}U_m^4.\nonumber
              \end{eqnarray*}
Using the definition of  $H_n(Z_i, Z_j)$, it is easy to see that $ \mathbb{E} \prod_{j=2}^l H_n(Z_1,Z_j)=0$ for any $l>1$. Therefore,
         \begin{eqnarray}
  \mathbb{E}Y_{nm}^4&=&\sum_{l_1=1}^{m-1}\sum_{l_2=1}^{m-1}\sum_{l_3=1}^{m-1}\sum_{l_4=1}^{m-1}  \mathbb{E} \prod_{j=1}^4H_n(Z_m,Z_{lj})\nonumber\\&=&\sum_{l=1}^{m-1} \mathbb{E} H_n(Z_m,Z_{l})^4 \nonumber \\
  &&+3 \sum_{1\leq l_1\neq l_2 \leq m-1} \mathbb{E} H_n(Z_m,Z_{l_1})^2H_n(Z_m,Z_{l_2})^2.
      \nonumber  \end{eqnarray}   Then using Holder's inequality and by setting $K_0=4$ in Lemma \ref{night}   we have \begin{eqnarray*}
&&\sum_{m=2}^{n-1}  \mathbb{E}Y_{nm}^4\\
&\leq&n^2\mathbb{E} H_n(Z_1,Z_2)^4+n^3 \mathbb{E} H_n(Z_1,Z_{2})^2H_n(Z_1,Z_{3})^2\\
&\leq&C\frac{p^4}{n^5},      \end{eqnarray*}
  $\sum_{m=1}^{n-1}  \mathbb{E}Y_{m}^4\leq C \frac{p^2}{n^3}$ and   $\sum_{m=1}^{n-1}  \mathbb{E}U_{m}^4 \leq C \frac{p^2\tau^4C_n^4}{n^7}$. Combining these with the lower bound of $s_{n-1}^2$ in Lemma \ref{son}, we have
\begin{eqnarray*}\lbl{use}
s_{n-1}^{-4}\sum_{m=2}^{n-1}   \mathbb{E}G_m^4 \to 0
  \end{eqnarray*}
  as $n\to\infty$.  This  implies (\ref{cat}).

We now prove (\ref{dog1}).
By the definition of $Y_{nm}$, observe that
                \begin{eqnarray*}
 \mathbb{E}_{m-1}G_m^2&=&  \mathbb{E}_{m-1}(Y_{nm}+Y_m+U_m)^2\\
 & =&    \mathbb{E}_{m-1}Y_{nm}^2+  2 \mathbb{E}_{m-1}Y_{nm}Y_m  + \mathbb{E}Y_m^2\nonumber\\
 &+& 2 \mathbb{E}_{m-1}Y_{nm}U_m+  2 \mathbb{E} Y_mU_m+ \mathbb{E}U_m^2
              \end{eqnarray*}
               and $s_{n-1}^2=\sum_{m=2}^{n-1} \mathbb{E}(Y^2_{nm})+(n-1) \mathbb{E}(Y^2_{1}+U^2_{1}+2Y_1U_1)$.
       Then $V_{n-1}^2-s_{n-1}^2=\sum_{m=2}^{n-1}[(\mathbb{E}_{m-1}Y_{nm}^2-\mathbb{E}Y_{nm}^2)+2\mathbb{E}_{m-1}Y_{nm}Y_m+2\mathbb{E}_{m-1}Y_{nm}U_m].$  Therefore it is easy to see that
       \begin{eqnarray*}\lbl{chen}
  &&  \mathbb{E}(V_{n-1}^2-s_{n-1}^2)^2 \\ &=&\mathbb{E}\{\sum_{m=2}^{n-1}(\mathbb{E}_{m-1}Y_{nm}^2 -\mathbb{E}Y_{nm}^2)\}^2 \nonumber \\
  &&+4\sum_{m=2}^{n-1}\sum_{k=2}^{n-1}\mathbb{E}\{\mathbb{E}_{m-1}(Y_{nm}Y_m)\mathbb{E}_{k-1}(Y_{nk}Y_k)\nonumber\\
  &&+\mathbb{E}_{m-1}(Y_{nm}U_m)\mathbb{E}_{k-1}(Y_{nk}U_k)\nonumber \\
  &&+\mathbb{E}_{m-1}(Y_{nm}Y_m)\mathbb{E}_{k-1}(Y_{nk}U_k)\nonumber\\
  &&+\mathbb{E}_{m-1}Y_{nm}^2\mathbb{E}_{k-1}(Y_{nk}Y_k) \mathbb{E}_{m-1}Y_{nm}^2\mathbb{E}_{k-1}(Y_{nk}U_k)\}\nonumber.
       \end{eqnarray*}
      Let $K_0=2$  in Lemma \ref{dress}, and by Lemma \ref{zhu} we see that
      \begin{eqnarray*}\lbl{yi}
       \mathbb{E}(V_{n-1}^2-s_{n-1}^2)^2 \leq C(  \frac{p^3}{n^5}+\frac{p(\tau C_n)^2}{n^5}+\frac{p^2 \tau C_n}{n^5}+\frac{p^2}{n^4})
.\end{eqnarray*} Using the lower bound $ \frac{1}{C}\max(\frac{p^2}{n^2}, \frac{p(\tau C_n)^2}{n^3})$  of $s_{n-1}^2$ from Lemma \ref{son},  it is easy to see
             $s_{n-1}^{-4}  \mathbb{E}(V_{n-1}^2-s_{n-1}^2)^2 \to 0$            as $n\to \infty$,     which implies (\ref{dog1}). Therefore we can conclude that $s_{n-1}^{-1}S_{n-1}$ converges to a standard normal distribution.

By Lemma \ref{son}, we know $s_{n-1}^2=\text{Var}_n(\tau)$ as given in (\ref{sam}).
By the  Slutsky's theorem, in order to derive the central limit theorem for $  \text{Var}_n(\tau)^{-\frac{1}{2}}(S_n^c(\tau)-R_c(\tau))$, we only need to show that $
  \text{Var}_n(\tau)^{-\frac{1}{2}} R_1\to 0$  and  $  \text{Var}_n(\tau)^{-\frac{1}{2}}R_2 \to 0$ in probability as $n\to \infty$.  Combining Lemmas \ref{son}, \ref{night} and \ref{night1},  it is easy to see that  when $n\to \infty$,
  \begin{eqnarray*}
  \text{Var}_n(\tau) ^{-2} \mathbb{E}R_1^4\, \text{ and}\,    \text{Var}_n(\tau)^{-2} \mathbb{E}R_2^4 \to 0.
   \end{eqnarray*} By Chebyshev's inequality,  we  get the desired result.
\end{IEEEproof}

\begin{IEEEproof}[Proof of Theorem \ref{babygirl}]
Since $C_n$ is a constant and $p_n/n\to y\in[1,\infty)$,   we know that  the lower bound of $s_{n-1}^2$ is $\frac{p^2}{Cn^2}$ from the proof of Lemma 9 (see (\ref{chi}) in the appendix). So using this new lower bound  in the proof of Theorem  \ref{girl},  keeping all other statements the same, we  get the desired result.
\end{IEEEproof}

\begin{IEEEproof}[Proof of Theorem \ref{zhen}]   Since $S_n^c(\tau)$ has the same distribution as $\text{SURE}_c(\tau)$, we only need to show that the probability bound for $S_n^c(\tau)-R_c(\tau)$. Following the same notation as in the proof of Theorem \ref{girl},  $S_n^c(\tau)-R_c(\tau)=S_{n-1}+R_1+R_2$, where $S_{n-1}=\sum_{m=1}^{n-1}G_m$. Set $h_n=\lambda\sqrt{\text{Var}_n(\tau)}/t_n$, $\mu_m=\mathbb{E}(G_mI(|G_m|\leq h_n)|\ml{F}_{m-1})$,
\begin{eqnarray}\lbl{nursery}\eta_m&=&G_m I(|G_m|\leq h_n)-\mu_m \nonumber \\
xi_m&=&G_m I(|G_m|> h_n)+\mu_m\end{eqnarray} for any $m\geq 1$.
Using the inequality $\mathbb{P}(U+V\geq u+v)\leq \mathbb{P}(U\geq u)+P(V\geq v)$, we obtain that
\begin{eqnarray*}\lbl{contract}
&&\mathbb{P}(  |S_n^c(\tau)-R_c(\tau)|\geq \lambda t_n \sqrt{\text{Var}_n(\tau)} )\nonumber\\
&\leq & \mathbb{P}( |S_{n-1}|\geq \frac{7}{8}\lambda t_n \sqrt{\text{Var}_n(\tau)} )\nonumber\\
&&+\mathbb{P}(|R_1+R_2|\geq\frac{1}{8}\lambda t_n\sqrt{\text{Var}_n(\tau)})\nonumber\\
&\leq &\mathbb{P}(|\sum_{m=1}^{n-1}\eta_m|>\frac{3}{4}\lambda t_n \sqrt{\text{Var}_n(\tau)})\nonumber\\
&&+\mathbb{P}(|\sum_{m=2}^{n-1}\xi_m|>\frac{1}{8}\lambda t_n \sqrt{\text{Var}_n(\tau)})\nonumber\\
&&+\mathbb{P}(|R_1+R_2|\geq \frac{1}{8}  \lambda t_n\sqrt{\text{Var}_n(\tau)})\nonumber\\
&:=& \ml{A}_n+\ml{B}_n+\ml{C}_n.\end{eqnarray*}

First, we find the bound of $\ml{A}_n$. Since $\{\eta_m, \ml{F}_m\}$ is a martingale difference sequence  with $\mathbb{E}(\eta_m |\ml{F}_{m-1})=0$, set $\sigma_m^2=\mathbb{E}(\eta_m^2|\ml{F}_{m-1})$   and $T_{n-1}^2=\sum_{m=1}^{n-1}\sigma_m^2$. Furthermore $|\eta_m|\leq 2 h_n$ for any $m\geq 1$. Then use Lemma \ref{DE} (also see De La Pe$\tilde{n}$a, 1999),  by letting $L_n=3\text{ Var} _n(\tau)=3s_{n-1}^2$,
\begin{eqnarray*} \ml{A}_n&\leq&\mathbb{P}(|\sum_{m=1}^{n-1}\eta_m|>\frac{3\lambda }{4}t_n \sqrt{\text{Var}_n(\tau)}, T_{n-1}^2\leq L_n)\\
&&+\mathbb{P}(T_{n-1}^2\geq L_n)   \nonumber\\
&\leq&2\exp\bigg(-\frac{(\frac{3\lambda }{4})^{2}t_n^2\text{Var}_n(\tau)}{2(L_n+h_n\frac{3\lambda }{2} t_n \sqrt{\text{Var}_n(\tau)})}\bigg)\\
&&+\mathbb{P}(T_{n-1}^2\geq L_n)\nonumber\\
&\leq&2\exp\big(-\frac{1}{16}t_n^2\big)+\mathbb{P}(\sum_{m=1}^{n-1}\mathbb{E}(G_m^2|\ml{F}_{m-1})\geq L_n).\end{eqnarray*}
 Notice that $\sum_{m=1}^{n-1}\mathbb{E}(G_m^2|\ml{F}_{m-1})\leq 2\sum_{m=2}^{n-1} \mathbb{E}_{m-1} Y_{nm}^2+2 \sum_{m=1}^{n-1}\mathbb{E}(Y_m+U_m)^2 $. Thus by Lemma \ref{son} and Lemma \ref{dress}, \begin{eqnarray*}
 &&\mathbb{P}(\sum_{m=1}^{n-1}\mathbb{E}(G_m^2|\ml{F}_{m-1})\geq L_n)\\
 &\leq&\mathbb{P}(2\sum_{m=2}^{n-1}(\mathbb{E}_{m-1}Y_{nm}^2-\mathbb{E}Y_{nm}^2)\geq s_{n-1}^2) \nonumber\\
 & \leq& \frac{2^{K_0}\mathbb{E}(\sum_{m=2}^{n-1}(\mathbb{E}_{m-1}Y_{nm}^2-\mathbb{E}Y_{nm}^2))^{K_0}}{ s_{n-1}^{2K_0}}\nonumber\\&\leq& M_{K_0}\big(\frac{1}{p^{K_0}}+\frac{1}{(np)^{K_0/2}}\big).\end{eqnarray*} Then
       \begin{eqnarray}\lbl{butterfly}  \ml{A}_n&\leq& 2\exp\big(-\frac{1}{16}t_n^2\big)+   M_{K_0}(np)^{-K_0/2}.\end{eqnarray}
       For computing the bound of $\ml{B}_n$, first recall the definition of $\xi_m$ and $\mu_m$ in (\ref{nursery}), we have,
       \begin{eqnarray*}\lbl{unlimited}  \ml{B}_n&\leq&  \frac{64\mathbb{E}(\sum_{m=1}^{n-1}\xi_m)^2}{\lambda^2t_n^2 \text{Var}_n(\tau)}\\
       &=& \frac{64\sum_{m=1}^{n-1}\mathbb{E} \xi_m^2}{\lambda^2t_n^2 \text{Var}_n(\tau)}\nonumber\\
       &\leq&\frac{64\sum_{m=1}^{n-1}\mathbb{E} G_m^2I (|G_m|\geq h_n)}{\lambda^2t_n^2 \text{Var}_n(\tau)}\\
       &\leq &\frac{64\sum_{m=1}^{n-1}\mathbb{E}( G_m)^{K_0} }{\lambda^2t_n^2 \text{Var}_n(\tau)h_n^{K_0-2}}\nonumber\\
       &\leq&  \frac{3^{K_0-1}64t_n^{K_0-4} }{\lambda^{K_0}\text{Var}_n(\tau)^{K_0/2}} \cdot\\
       &&\left\{\sum_{m=2}^{n-1} \mathbb{E} Y_{nm}^{K_0}+(n-1)[\mathbb{E}U_1^{K_0}+\mathbb{E}Y_1^{K_0}]\right\},
       \end{eqnarray*}  for any even number $K_0\geq 4$. By Moment inequality for independent variables and apply Lemma \ref{night} , we have that
 \begin{eqnarray*}
 &&\mathbb{E}Y_{nm}^{K_0}\\
 &=&\mathbb{E}\big\{\mathbb{E}\{[\sum_{l=1}^mH_n(Z_m, Z_l)]^k|Z_m\}\big\}\\
 &\leq &C(m-1)^{\frac{k_0}{2}} \mathbb{E}H_n(Z_1,Z_2)^{K_0}\leq C\frac{p^{K_0}}{n^{3K_0/2-1}},\end{eqnarray*}
  $   \mathbb{E}Y_{1}^{K_0}=C\frac{p^{K_0/2}}{n^{K_0 }}$ and   $  \mathbb{E}U_{1}^{K_0 }=C\frac{(C_n\tau)^{K_0 }p^{K_0/2}}{n^{2K_0}}$. By Lemma \ref{son}, \begin{eqnarray}\lbl{unlimited} \ml{ B}_n&\leq&  M_{K_0}\frac{ t_n^{K_0-4 }}{\lambda ^{K_0 }n^{K_0/2-1}}.
              \end{eqnarray}
Finally, we  calculate the bound of $\ml{C}_n.$ By Lemmas \ref{son} and \ref{night1}. it is easy to show that
\begin{eqnarray}\lbl{visa}   \ml{C}_n & \leq & \mathbb{P}(|R_1+R_2|\geq \frac{\lambda}{8} t_n  \sqrt{\text{Var}_n(\tau)} ) \nonumber \\
&\leq& \frac{16^{K_0}(\mathbb{E}R_1^{K_0}+\mathbb{E}R_2^{K_0})}{\lambda^{K_0}t_n^{K_0}\text{Var}_n(\tau)^{K_0/2}}\nonumber\\&\leq& M_{K_0}\big(\frac{1}{\lambda^{K_0}t_n^{K_0}n^{K_0/2}}+\frac{C_n^{K_0}}{\lambda^{K_0}t_n^{K_0}n^{  K_0}}\big).\end{eqnarray}
 Now combine (\ref{butterfly}), (\ref{unlimited}) and (\ref{visa}), $\ml{A}_n+\ml{B}_n+\ml{C}_n\leq 2 \exp\big(-\frac{1}{16}t_n^2\big)+    M_{K_0}\big\{(np)^{-K_0/2} + (\lambda^2 t_n^2n)^{-K_0/2}+  \frac{ t_n^{K_0-4 }}{\lambda ^{K_0 }n^{K_0/2-1}}+\frac{C_n^{K_0}}{(\lambda t_n)^{K_0}n^{ K_0}} )\big\}$. We arrive at the conclusion.
    \end{IEEEproof}

\begin{IEEEproof}[Proof of Theorem \ref{cooked}] For any $\Sigma\in \cal{F}_{\alpha}$,  by the definitions of $\hat{\tau}_{n}$ and $\tau_0$, we know that $\text{SURE}(\hat{\tau}_{n})\leq \text{SURE}(\tau_0)$. Then by (\ref{gen0}), it is easy to see that
\begin{eqnarray*}
&&\|\hat{\Sigma}^{(\hat{\tau}_n)}-\tilde{\Sigma}^{s}\|_{F}^2 \\
&\leq & \|\hat{\Sigma}^{(\tau_0)}-\tilde{\Sigma}^{s}\|_{F}^2+2\frac{n-1}{n}\sum_{i, j}\omega_{ij}^{(\tau_0)}\widehat{\text{var}}(\tilde{\sigma}_{ij}^s)\\
&&-2\frac{n-1}{n}\sum_{i, j}\omega_{ij}^{(\hat{\tau}_n)}\widehat{\text{var}}(\tilde{\sigma}_{ij}^s).
\end{eqnarray*}
So we can conclude that for any $\Sigma\in \cal{F}_{\alpha}$,
  \begin{eqnarray}\lbl{nuli}
  &&\mathbb{E}\|\hat{\Sigma}^{(\hat{\tau}_n)}-\Sigma\|_{F}^2 \nonumber \\
  &\leq& R(\tau_0) +2\frac{n-1}{n}\mathbb{E}\{\sum_{i, j}\omega_{ij}^{(\hat{\tau}_n)}\sigma_{ij}(\tilde{\sigma}_{ij}^{s}-\sigma_{ij})\}\nonumber\\&+&2\frac{n-1}{n}\mathbb{E}\sum_{i, j}\omega_{ij}^{(\hat{\tau}_n)}\big\{(\tilde{\sigma}_{ij}^{s}-\sigma_{ij})^2-\frac{\sigma_{ii}\sigma_{jj}+\sigma_{ij}^2}{n-1}\big\}\nonumber\\&-&2\frac{n-1}{n}\mathbb{E}\sum_{i, j}\omega_{ij}^{(\hat{\tau}_n)}\big\{\widehat{\text{var}}(\tilde{\sigma}_{ij}^2)-\frac{\sigma_{ii}\sigma_{jj}+\sigma_{ij}^2}{n-1}\big\}.\end{eqnarray}
Since   $\mathbb{E}\{\sum_{i, j}\omega_{ij}^{(k)}\sigma_{ij}(\tilde{\sigma}_{ij}^{s}-\sigma_{ij})\}^2\leq \frac{Cp}{n-1}$ for any $k\geq 1$,  we have
\begin{eqnarray*}\lbl{Hun}
&&|\mathbb{E}\{\sum_{i, j}\omega_{ij}^{(\hat{\tau}_n)}\sigma_{ij}(\tilde{\sigma}_{ij}^{s}-\sigma_{ij})\}|\nonumber \\
&=&\sum_{k=1}^p\mathbb{E}|\{\sum_{i, j}\omega_{ij}^{(k)}\sigma_{ij}(\tilde{\sigma}_{ij}^{s}-\sigma_{ij})I(\hat{\tau}_n=k)\}|\nonumber \\
&\leq&C\sqrt{\frac{p}{n}}\sum_{k=1}^p\{\mathbb{ P}(\hat{\tau}_n= k)\}^{\frac{1}{2}}.
\end{eqnarray*}
 To deal with the third term in (\ref{nuli}),  we have
  \begin{eqnarray*}&&\big|\mathbb{E}\sum_{i, j}\omega_{ij}^{(\hat{\tau}_n)}\big\{(\tilde{\sigma}_{ij}^{s}-\sigma_{ij})^2-\frac{\sigma_{ii}\sigma_{jj}+\sigma_{ij}^2}{n-1} \big\}\big|\nonumber\\
  &\leq&\sqrt{2} \sum_{k=1}^p\{ \mathbb{P}(\hat{\tau}_n= k)\}^{\frac{1}{2}}\big[\\
  &&\frac{1}{(n-1)^3}\mathbb{E}\{\sum_{i, j}\omega_{ij}^{(k)}(Z_{1i}Z_{1j}-\sigma_{ij})^2-\sigma_{ii}\sigma_{jj}-\sigma_{ij}^2\}^2\\
  &&+\frac{2(n-1)(n-2)}{(n-1)^4} \times \\
  &&\sum_{1\leq i, j, s, t \leq p}\omega_{ij}^{(k)}\omega_{st}^{(k)}\{\mathbb{E}(Z_{1i}Z_{1j}-\sigma_{ij})(Z_{1s}Z_{1t}-\sigma_{st})\}^2\\
  &&\big]^{\frac{1}{2}}\end{eqnarray*}
Furthermore, we know that \begin{eqnarray*}
 &&\mathbb{E}[\sum_{i, j}\omega_{ij}^{(k)}\{(Z_{1i}Z_{1j}-\sigma_{ij})^2-\sigma_{ii}\sigma_{jj}-\sigma_{ij}^2\}]^2\\
 & \leq &Cp^2k^2;\\
&& \sum_{1\leq i, j, s, t \leq p}\omega_{ij}^{(k)}\omega_{st}^{(k)}\{\mathbb{E}(Z_{1i}Z_{1j}-\sigma_{ij})(Z_{1s}Z_{1t}-\sigma_{st})\}^2\\
 &\leq& Cpk^2.
 \end{eqnarray*}
 So $\big|2\frac{n-1}{n}\mathbb{E}\sum_{i, j}\omega_{ij}^{(\hat{\tau}_n)}\{(\tilde{\sigma}_{ij}^{s}-\sigma_{ij})^2-\frac{\sigma_{ii}\sigma_{jj}+\sigma_{ij}^2}{n-1}\}\big|\leq C\frac{pk}{n^{3/2}}$.
 Then by the similar arguments, we know that $\big|2\frac{n-1}{n}\mathbb{E}\sum_{i, j}\omega_{ij}^{(\hat{\tau}_n)}\{\widehat{\text{var}}(\tilde{\sigma}_{ij}^2)-\frac{\sigma_{ii}\sigma_{jj}+\sigma_{ij}^2}{n-1}\}\big|\leq C\frac{pk}{n^2}$.
 So
 \begin{eqnarray}\lbl{wet}
&& \mathbb{E}\|\hat{\Sigma}^{(\hat{\tau}_n)}-\Sigma\|_{F}^2 \\
&\leq& R(\tau_0) +C\sqrt{\frac{p}{n}}\sum_{k=1}^p\{ \mathbb{P}(\hat{\tau}_n= k)\}^{\frac{1}{2}} \nonumber \\
&&+C\frac{p}{n^{\frac{3}{2}}}\sum_{k=1}^pk\{\mathbb{ P}(\hat{\tau}_n= k)\}^{\frac{1}{2}}.\nonumber
 \end{eqnarray}
 By Lemma \ref{Caizz}, we know that $$\sup\limits_{\Sigma\in\cal{F}_{\alpha}}R(\tau_0) \leq  \sup\limits_{\Sigma\in \mathcal{F}_{\alpha}}\mathbb{E}\|\Sigma^{({\tau}_n)}-\Sigma\|_{F}^2=I_n,$$ where $\tau_n=n^{1/2(\alpha+1)}$ and $I_n=C_0p n^{-(2\alpha+1)/2(\alpha+1)}$.

 We only need to find the upper bound of the rest terms in (\ref{wet}).
 Without loss of generality, we assume $\sigma_{ii}=1$ for all $1\leq i\leq p.$  Now choose $k_n=4C_1 n^{1/2(\alpha+1)}$, where $C_1=\max\{C_0, 1\}$. It is easy to see that $R(k)\geq \sum_{|i-j|\leq \lfloor \frac{k}{2}\rfloor}\frac{n-1}{n^2} \sigma_{ii}\sigma_{jj}=\frac{p(n-1)k}{n^2}\geq 3I_n$ for all $ k_n\leq k\leq \frac{p}{3}$ if $n\geq 4$.  From now on, we assume  $n$ is large enough. First, for any $k_n\leq k \leq \frac{p}{3}$,
 \begin{eqnarray*}
 &&\mathbb{P}(\hat{\tau}_{n}=k) \\
  &\leq &\mathbb{P}(|\text{SURE}(k)-R(k)|\geq \frac{R(k)}{3})\\
  &&+\mathbb{P}(|\text{SURE}(\tau_0)-R(\tau_0)|\geq \frac{R(k)}{3}).\end{eqnarray*}
Since $R(\tau_0)\leq I_n$, we know that  $\tau_0\leq C_0 n^{\frac{1}{2(\alpha+1)}}$. So by Lemma \ref{son}, there exists $C_2$ such that $\text{Var}_n(\tau_0)\leq C_2\frac{(n-1)^2p^2}{n^4}$,$\text{Var}_n(k)\leq C_2\frac{(n-1)^2p^2}{n^4}$ for $k_n\leq k\leq\frac{1}{3} \sqrt{np}$ and  $\text{Var}_n(k)\leq C_2\frac{(n-1)^2pk^2}{n^5}$ for $ \frac{1}{3}\sqrt{np}\leq k\leq \frac{p}{3}$. Now choose $K_0=8$ and apply Theorem \ref{zhen} with different $t_n$ and $\lambda$   to bound $\mathbb{P}(|\text{SURE}(k)-R(k)|\geq \frac{R(k)}{3})$  for different $k $.  For $k_n\leq k\leq \frac{n}{3}$, set $t_n=\sqrt{\frac{k_n}{ C_2}} $ and $\lambda_n=\frac{k}{3\sqrt{k_n}}\geq 1$, there exists a constant $C$ such  that
\begin{eqnarray}\lbl{tian}
&&\mathbb{P}(|\text{SURE}(k)-R(k)|\geq \frac{R(k)}{3})\\
&\leq &2\exp\big(-\frac{k_n}{16C_2}  \big)+C\frac{ (k_n)^{6}}{k^{8}n^{3}}.\nonumber
\end{eqnarray}
For $\frac{n}{3} \leq k\leq \frac{\sqrt{np}}{3} $,  choose $t_n=\frac{(\log p)^{5/8}}{\sqrt{C_2}} $ and $\lambda_n=\frac{k}{3(\log p)^{5/8}}\geq 1$,
there exists $C$ such that
\begin{eqnarray}\lbl{di}
&&\mathbb{P}(|\text{SURE}(k)-R(k)|\geq \frac{R(k)}{3}) \\
&\leq & 2\exp\big(-\frac{(\log p)^{5/4}}{16C_2}  \big)+ C\big\{(np)^{-4}+\frac{(\log p)^{15/2}}{ k^{8}n^{3}}\big\}.  \nonumber
\end{eqnarray}
For $  (np)^{\frac{1}{2}}\leq k\leq \frac{p}{3}$,    choose $t_n=\frac{\log p}{\sqrt{C_2}}$ and $\lambda=\frac{(np)^{\frac{1}{2}}}{\log p}$,  there exists a constant $C$
\begin{eqnarray}\lbl{jun}
&&\mathbb{P}(|\text{SURE}(k)-R(k)|\geq \frac{R(k)}{3})\\
&\leq& 2\exp\big(-\frac{(\log p)^2}{16C_2}  \big)+C\big \{(np)^{-4}+\frac{(\log p )^{12}}{n^{7}p^{4}}\big\} \nonumber.\end{eqnarray}
The same results from (\ref{tian}), (\ref{di}) and (\ref{jun}) are also true for $\mathbb{P}(|\text{SURE}(\tau_0)-R(\tau_0)|\geq \frac{R(k)}{3})$.
For any $\tau$, by applying Lemmas \ref{night} and \ref{night1},  we have that, for any even number $K_0$,
\begin{eqnarray}\lbl{backpain}
&&\frac{n^{K_0}\mathbb{E}(S_n^c(\tau)-R(\tau))^{K_0}}{ p^{2K_0} }\nonumber \\
&\leq&  \frac{3^{K_0-1} L_{K_0}n^{K_0}[n^{\frac{K_0-2}{2}} \sum_{m=1}^{n}\mathbb{E}G_m^{K_0}+\mathbb{E}R_1^{K_0}+\mathbb{E}R_2^{K_0}]}{p^{2K_0} }\nonumber\\
&\leq&C(np)^{-K_0/2}. \quad  \end{eqnarray}

For $k\geq \frac{p}{3}$, we know that $R(k)\geq \frac{p^2}{3n}$. Then for $K_0= 8$,  by Markov inequality, $ \mathbb{P}(\hat{\tau}_{n}=k)\leq C(np)^{-4}$. Now we have
     \begin{eqnarray}
\{\mathbb{P}(\hat{\tau}_n=k)\}^{\frac{1}{2}}\leq\left\{ \begin{array}{ll}2e^{-\frac{k_n}{32C_2}}+C\frac{ k_n^{3}}{k^{4}n^{3/2}},   \mbox{ $k_n\leq k< \frac{n}{3}$};\nonumber\\
      2e^{-\frac{(\log p)^{\frac{5}{4}}}{32C_2}}+\frac{C(\log p)^{15/4}}{k^{4}n^{3/2}} +\frac{C}{n^2p^2} ,    \mbox{ $\frac{n}{3} \leq k< \frac{ \sqrt{np}}{3} $};\nonumber\\
  2e^{-\frac{(\log p)^2}{32C_2}}+\frac{C1}{n^2p^2}+\frac{C (\log p)^6}{p^{2}n^{7/2}} , \mbox{ $\frac{ \sqrt{np}}{3} \leq  k< \frac{p}{3}$};\nonumber\\
 \frac{C}{(np)^{2}}, \mbox{$k
 \geq \frac{p}{3}$}.\end{array} \right.
         \end{eqnarray}
         So when $n$ is large enough, we have that
         \begin{eqnarray}\lbl{miao} \sum_{k=k_n}^{p}\frac{k}{\sqrt{n}}\{\mathbb{P}(\hat{\tau}_n=k)\}^{\frac{1}{2}}\leq Cn^{-\frac{1}{
       4}}\, \sum_{k=k_n}^{p}\{\mathbb{P}(\hat{\tau}_n=k)\}^{\frac{1}{2}}\leq C n^{-\frac{3}{
       4}}. \end{eqnarray}
 It is easy to check that $\sum_{k=1}^{k_n}\{\mathbb{P}(\hat{\tau}_n= k)\}^{\frac{1}{2}}\leq k_n^{1/2}.
$  So  \begin{eqnarray}\lbl{wet2}
  \big[\sqrt{\frac{p}{n}}\sum_{k=1}^p\{\mathbb{P}(\hat{\tau}_n= k)\}^{\frac{1}{2}}+C\frac{p}{n^{\frac{3}{2}}}\sum_{k=1}^pk\{\mathbb{P}(\hat{\tau}_n= k)\}^{\frac{1}{2}}\big]/I_n\to 0
 \end{eqnarray} as $n\to \infty$.
Then we see that $\mathbb{E}\|\hat{\Sigma}^{(\hat{\tau}_n)}-\Sigma\|_{F}^2\asymp pn^{-(2\alpha+1)/2(\alpha+1)}$ for any $\Sigma\in \mathcal{F}_{\alpha}$.
\end{IEEEproof}

\begin{IEEEproof}[Proof of Theorem \ref{daughter}]
  For any $|h|\geq 1$, it is easy to check that
\begin{eqnarray}\lbl{home}
&&\mathbb{P}(\hat{\tau}_n=\tau_0+h)\nonumber\\
&\leq& \mathbb{P}(\text{SURE}(\tau_0+h)< \text{SURE}(\tau_0))\nonumber\\
&\leq&  \mathbb{P}(\text{SURE}(\tau_0+h)-R(\tau_0+h)< -\frac{R(\tau_0+h)-R(\tau_0)}{2})\nonumber\\&&+\mathbb{P}(\text{SURE}(\tau_0)-R(\tau_0)>\frac{R(\tau_0+h)-R(\tau_0)}{2})
\nonumber\\&:=&A_n+B_n.
 \end{eqnarray}
Since $\tau_0$ is the unique minimizer of $R(\tau)$, we know that $R(\tau_0\pm1)-R(\tau_0)>0$. Then, \begin{eqnarray}\lbl{cough}\sum_{|i-j|=\tau_0}\sigma_{ij}^{2}< \sum_{|i-j|=\tau_0}\frac{1}{n}\sigma_{ii}\sigma_{jj}<\sum_{|i-j|=\tau_0-1}\sigma_{ij}^{2}. \end{eqnarray}
   By (C.3) conditions, we know that $  |\sigma_{ij}|\leq M_1 (\tau)^{-(\alpha+1)}$ when $|i-j|=\tau$ for all $1\leq i\leq p$. So there exists a constant C such that $\tau_0\leq C n^{\frac{1}{2(\alpha+1)}} \ll  (np)^{\frac{1}{2}}$.
Furthermore $\sum_{|i-j|=k}\frac{(n-1)}{n}=\frac{2p(n-1)}{n}$ as long as $k<\frac{p}{2}.$ For any $h\leq \frac{p}{3}$,
if $\sum_{|i-j|=\tau_0}\sigma_{ij}^2\asymp \frac{p  }{n}$, then
 \begin{eqnarray}\lbl{tear1}R(\tau_0+h)-R(\tau_0)\geq (|h|-1)(1-\frac{1}{\gamma})\sum_{|i-j|=\tau_0}\sigma_{ij}^2;\end{eqnarray}
 if  $\sum_{|i-j|=\tau_0}\frac{n+1-\log n}{n}\sigma_{ij}^2\ll \frac{p }{n}$, then   \begin{eqnarray}\lbl{dcare1}R(\tau_0+h)-R(\tau_0)\geq  (|h|-1)(\gamma-1)  \frac{2(n-1)p}{n^2}. \end{eqnarray}
 Then there exists $\delta>0$ such that
  $\frac{R(\tau_0+h)-R(\tau_0)}{2}\geq |h|\delta \frac{p}{n}$ for $ 2\leq |h |\leq \frac{p}{3}$.
By Lemma \ref{son}, there exists a constant C such that $\text{Var}_n(\tau_0)\leq C\frac{p^2}{n^2} $ and $\text{Var}_n(\tau_0+h)\leq C\frac{p^2}{n^2}$ for $|\tau_0+h|\leq \frac{\sqrt{np}}{3} $.  From now on, we assume  $n$ is large enough and apply Theorem \ref{zhen} to bound $A_n$ and $B_n$ with different $t_n$ and $\lambda$ for different $h$. We further choose $K_0=8$.  For $\log n \leq |h|\leq \log p$ and $\tau_0+h\geq 1$,   set $t_n=\frac{\delta\log n}{\sqrt{C}}$ and $\lambda= \frac{|h|}{\log n}$,     there exist a constant $C$ such that the upper bound for $A_n$ and $B_n$  is $  2\exp\big(-\frac{\delta^2}{16C}(\log n)^2\big)+    C \frac{  ( \log n)^{12 }}{ |h| ^{8 }n^{3}}.$
For  $\log p\leq |h|$ and $|h+\tau_0|\leq  \frac{\sqrt{np}}{3}$ and $\tau_0+h\geq 1$,  choose $t_n=\frac{\delta(\log p)^{ 5/8}}{\sqrt{C}}$ and $\lambda= \frac{|h|}{(\log p)^{5/8}}$,   the bound is \begin{eqnarray*}  2 \exp\big(-\frac{\delta^2}{16C}(\log p)^{5/4}\big)+    C\big(\frac{ ( \log p)^{15/2 }}{|h| ^{8 }n^{3}}+(np)^{-4}\big).\end{eqnarray*}   For $\frac{\sqrt{np}}{3}\leq \tau_0+h\leq \frac{p}{3}$,  by choosing $t_n=\frac{\delta\log p}{\sqrt{C}}$ and $\lambda=\frac{h}{\log p}$, we have that
\begin{eqnarray*}  B_n \leq 2 \exp\big(-\frac{\delta^2}{16C}(\log p)^{2}\big)+    C\big( (np)^{-4}+  \frac{   ( \log p)^{12 }}{h ^{8 }n^{3}}\big).\end{eqnarray*}
By Lemma \ref{son},
therefore there exist $C_1>1$ such that
  $ \frac{ p^2}{n^2 \text{Var}_n(\tau_0+h )}\geq\frac{np}{C_1(\tau_0+h)^2}$ for $\frac{\sqrt{np}}{3}\leq \tau_0+h\leq \frac{p}{3}$. Since $n$ is large enough, we have that $\frac{h}{\tau_0+h}\geq \frac{1}{\sqrt{2}}$. Then $A_n=\mathbb{P}(|\text{SURE}(\tau_0+h)-R(\tau_0+h)|\geq\frac{\delta\sqrt{np}{\text{Var}_n(\tau_0+h)}^{\frac{1}{2}}}{\sqrt{2C_1}}).$
Applying Theorems \ref{zhen} to $A_n$ with $K_0=8$, $t_n=\frac{\delta}{\sqrt{2C}}\log p$ and $\lambda=\frac{\sqrt{np}}{\log p}$,  we get \begin{eqnarray*} A_n \leq 2 \exp\big(-\frac{\delta^2}{32C}(\log p)^2\big)+    C\big( (np)^{-4} +  \frac{(\log p)^{12}}{ n^{7 }p^{4}}\big).\end{eqnarray*} When $h+\tau_0\geq \frac{p}{3}$, then $\frac{R(\tau_0+h)-R(\tau_0)}{2}\geq\frac{1}{C} \frac{p^2 }{n}$.
By (\ref{backpain}) and   let $K_0=8$ using Markov inequality, we get   $A_n +B_n\leq C(np)^{-4}$.
 , we can conclude that \begin{eqnarray}\lbl{tea1}
\sum_{|k- \tau_0|> \log n} \mathbb{P}(\hat{\tau}_n=k) &\leq& Cn^{-5/2},\end{eqnarray} when $n $ is large enough.     So $\sum_{n=1}^{\infty}\mathbb{P}(|\hat{\tau}_n- \tau_0|>\log n)<\infty$,
  by Borel--Cantelli Lemma, we know that $|\hat{\tau}_n -\tau_0|\leq \log n$ almost surely as $n\to \infty$.

By definition, $R(\tau_0)=\inf\limits_{\tau}\mathbb{E}\|\hat{\Sigma}^{(\tau)}-\Sigma\|_{F}^2$, $R(\tau_0)\leq \mathbb{E}\|\hat{\Sigma}^{(\hat{\tau}_n)}-\Sigma\|_{F}^2$.  By (\ref{vapor}), we have  that $R(\tau_0)\geq C\frac{p\tau_0}{n}$. So combining these with (\ref{wet}), it is easy to see that
  \begin{eqnarray*}\lbl{shoot}
  &&1 \leq\frac{ \mathbb{E}\|\hat{\Sigma}^{(\hat{\tau}_n)}-\Sigma\|_{F}^2 }{R(\tau_0)}\nonumber \\
  &\leq& 1 +C\sqrt{\frac{n}{p}}\frac{1}{\tau_0}\sum_{k=1}^{p} \{\mathbb{ P}(\hat{\tau}_n= k)\}^{\frac{1}{2}}\\
  &&+\frac{1}{\tau_0}\sum_{k=1}^{p}\frac{k}{\sqrt{n}}\{\mathbb{ P}(\hat{\tau}_n= k)\}^{\frac{1}{2}}.  \nonumber
  \end{eqnarray*}
From (\ref{miao}), we only need to show that $C\sqrt{\frac{n}{p}}\frac{1}{\tau_0}\sum_{k=1}^{k_n} \{\mathbb{P}(\hat{\tau}_n= k)\}^{\frac{1}{2}}+\frac{1}{\tau_0}\sum_{k=1}^{k_n}\frac{k}{\sqrt{n}}\{\mathbb{P}(\hat{\tau}_n= k)\}^{\frac{1}{2}}\to 0$ as $n\to \infty$. Following the similar proof of (\ref{tea1}) and choosing  $K_0=8$, we  can get that $\sum\limits_{|k-\tau_0|\geq \log n,k\leq k_n}\{\mathbb{P}(\hat{\tau}_n= k)\}^{\frac{1}{2}}\leq Cn^{-\frac{3}{4}}$. Furthermore, $\sum\limits_{|k-\tau_0|\leq \log n}\{\mathbb{P}(\hat{\tau}_n=k)\}^{\frac{1}{2}}\leq C\sqrt{\log n}.$  So $\frac{ \mathbb{E}\|\hat{\Sigma}^{(\hat{\tau}_n)}-\Sigma\|_{F}^2 }{R(\tau_0)}\leq 1+C\sqrt{\frac{n\log n}{p\tau_0^2}}+C\frac{(\tau_0+\log n)\log n}{\tau_0\sqrt{n}}$. Then if $p\gg n\log n$, we get the desired conclusion.
\end{IEEEproof}

\begin{IEEEproof}[Proof of Theorem \ref{Dad}]
By the definition of $R(\tau)$, we have that $R(\tau+
1)-R(\tau)=\sum_{|i-j|=\tau}\frac{n-1}{n^2}\sigma_{ii}\sigma_{jj}-(1-\frac{1}{n})\sigma_{ij}^2$. Then it is easy to see that $k_0$ is the unique minimizer of $R(\tau)$, because
  the off-diagonal entries of the covariance matrix equal zero when $|i-j|\geq k_0$ and the rest are bounded away from 0. For the simplicity of the proof, we assume $\sigma_{ii}=1$. Since $k_0$ is a constant, then  exact banded matrix $\Sigma\in \cal{G}_{\alpha}$ with $M_1=k_0^2$ and $\alpha<1.$ By Lemma \ref{son}, we know that $\text{Var}_n(k_0)\leq C\frac{p^2}{n^2}$. So  For $h\geq 1$, there exists $\delta>0$ such that $R(k_0+h)-R(k_0)\geq \sum_{k_0\leq |i-j|< k_0+h}\frac{n-1}{n^2}\sigma_{ii}\sigma_{jj}\geq 2 h\delta \sqrt{\text{Var}_n(k_0)}$ and $ R(k_0-h)-R(k_0)\geq \sum_{k_0-h\leq |i-j|\leq k_0-1}(\frac{n-1}{n})\sigma_{ij}^2-\frac{n-1}{n^2}\sigma_{ii}\sigma_{jj}\gg 2\delta h \log n\sqrt{\text{Var}_n(k_0)}$. Following the proof of Theorem \ref{daughter},   we know that $|\hat{\tau}_n-\tau_0|\leq \log n$ almost surely as $n\to \infty.$ By Borel-Cantelli Lemma, we only need to show that $\sum_{n=1}^{\infty}\sum_{k=1}^{k_0-1}P(\hat{\tau}_n=k)<\infty.$ By (\ref{home}), we only need to bound $A_n$ and $B_n$ for all $k\leq k_0-1$. Applying Theorem \ref{zhen} with $K_0=6$, $t_n=\delta \log n$ and $\lambda=h$, we find the upper bound $  2\exp\big(-\frac{\delta^2}{16C}(\log n)^2\big)+    C \frac{   ( \log n)^{2 }}{h^{6 }n^{3}}.$ Hence we get the desired conclusion.
\end{IEEEproof}

\begin{IEEEproof}[Proof of Theorem \ref{bake}]
Similar to the proof of Theorem \ref{daughter}, we want to find the upper bound of $A_n$ and $B_n$ in (\ref{home}). By the similar arguments as in (\ref{tear1}) and (\ref{dcare1}), we have that:
 if $\sum_{|i-j|=\tau_0}\sigma_{ij}^2\asymp \frac{p \log n}{n}$,
 \begin{eqnarray*}\lbl{tear}
 &&R_c(\tau_0^c+h)-R_c(\tau_0^c)\\
 &\geq& (|h|-1)(1-\frac{1}{\gamma})\sum_{|i-j|=\tau_0^c}\sigma_{ij}^2;\end{eqnarray*}
 if  $\sum_{|i-j|=\tau_0}\frac{n+1-\log n}{n}\sigma_{ij}^2\ll \frac{p\log n}{n}$,
 \begin{eqnarray*}\lbl{dcare}
 &&R_c(\tau_0^c+h)-R_c(\tau_0^c) \\
 &\geq &2(|h|-1)(\gamma-1)  \frac{(n-1)p+p(n(\log n -2))}{n^2}.
 \end{eqnarray*}
  Then there exists $\delta>0$ such that
  $\frac{R_c(\tau_0^c+h)-R_c(\tau_0^c)}{2}\geq |h|\delta \frac{p\log n}{n}$ for $ 2\leq |h|\leq \frac{p}{3}$.     Since  $\tau_0^c\leq C(\frac{n}{\log n})^{\frac{1}{2(\alpha+1)}} \ll \frac{(np)^{\frac{1}{2}}}{\log n}$,
    using Lemma \ref{son},  there exists a constant C such that $\text{Var}_n(\tau_0^c)\leq C\frac{p^2}{n^2} $ and $\text{Var}_n(\tau_0^c+h)\leq C\frac{p^2}{n^2}$ for $|\tau_0^c+h|\leq \frac{\sqrt{np}}{\log n} $.  From now on, we assume $n$ is large enough. Choose $K_0=8$ and apply Theorem \ref{zhen} to bound $A_n$ and $B_n$ with different $t_n$ and $\lambda$ for different $h$. For  $2 \leq |h|\leq \log p$ and $\tau_0^c+h\geq 1$,   with  $t_n=\frac{\delta\log n}{\sqrt{C}}$ and $\lambda= |h|$,   it is easy to see that  the upper bound for $A_n$ and $B_n$  is \begin{eqnarray*}  2\exp\big(-\frac{\delta^2}{16C}(\log n)^2\big)+    C   \frac{    ( \log n)^{4 }}{  |h| ^{8 }n^{3}} .\end{eqnarray*}
For  $\log p\leq |h|$ and $h+\tau_0^c\leq  \frac{\sqrt{np}}{\log n}$,  set  $t_n=\frac{\delta(\log p)^{5/8}\log n}{\sqrt{C}}$ and $\lambda= \frac{|h|}{(\log p)^{5/8}}$,   the bound is \begin{eqnarray*}  2 \exp\big(-\frac{(\delta\log n)^2}{16C}(\log p)^{\frac{5}{4}}\big)+    C\big((np)^{-4} +  \frac{    ( \log p)^{15/2 }}{ |h| ^{8 }n^{3}}\big).\end{eqnarray*} Furthermore for $\frac{\sqrt{np}}{\log n}\leq \tau_0^c+h\leq \frac{p}{3}$,  by choosing $t_n=\frac{\delta\log p}{\sqrt{C}}$ and $\lambda=\frac{|h|\log n}{\log p}$, we know
\begin{eqnarray*}  B_n \leq 2 \exp\big(-\frac{\delta^2}{16C}(\log p)^{2}\big)+   2 M_{K_0}(np^{-4} +  \frac{   ( \log p)^{12}}{  (|h|\log n) ^{8 }n^{3}}).\end{eqnarray*}
Using Lemma \ref{son},
therefore there exist $C>1$ such that
  $ \frac{ p^2}{n^2 \text{Var}_n(\tau_0^c+h )}\geq\frac{np}{C(\tau_0^c+h)^2(\log n)^2}$  for $\frac{\sqrt{np}}{3}\leq \tau_0^c+h\leq \frac{p}{3}$. Since $n$ is large enough, we have that $\frac{h}{\tau_0^c+h}\geq \frac{1}{\sqrt{2}}$  then
 $A_n= \mathbb{P}(|\text{SURE}_c(\tau_0^c+h)-R_c(\tau_0^c+h)|\geq \frac{\delta\sqrt{np}{\text{Var}_n(\tau_0^c+h)}^{\frac{1}{2}}}{\sqrt{2C}}).$
 Now  applying Theorem \ref{zhen} to $A_n$ with choosing $t_n=\frac{\delta}{\sqrt{2C}}\log p$ and $\lambda=\frac{\sqrt{np}}{\log p}$,  we get
 \begin{eqnarray*} A_n \leq2 \exp\big(-\frac{\delta^2}{32C}(\log p)^2\big)+    C\big( (np)^{-4} +  \frac{    ( \log p)^{12 }}{  n^{7 }p^{4}}\big).\end{eqnarray*} When $h+\tau_0\geq \frac{p}{3}$, then $\frac{R(\tau_0^c+h)-R(\tau_0^c)}{2}\geq\frac{1}{C} \frac{p^2 }{n}$.
Let $K_0=8$, by (\ref{backpain}), we get   $A_n +B_n\leq C(np)^{-4}$.
   We can conclude that \begin{eqnarray}\lbl{teaparty}
\sum_{|k- \tau_0^c| > 1} \mathbb{P}(\hat{\tau}_n^c=k) &\leq& Cn^{-5/2},\end{eqnarray} when $n $ is large enough.     So
  by Borel--Cantelli Lemma, we know that   $|\hat{\tau}_n^{c}-\tau_0^c|\leq 1$ almost surely. With further condition $|R_c(\tau_0^c\pm 1)-R_c(\tau_0^c)|\geq 2\delta \log n\sqrt{\text{Var}_n(\tau_0^c)}$, by the similar argument, we conclude that $\hat{\tau}_n^{c}$ converges to $\tau_0^c$ almost surely.
\end{IEEEproof}

\begin{IEEEproof}[Proof of Theorem \ref{breaking}] With Theorem \ref{bake} in hand, the proof is similar to the proof of Theorem~\ref{Dad}. For the sake of space we omit the proof here.
\end{IEEEproof}

\subsection{Appendix}
In this appendix we give the proofs of some technical Lemmas. In what follows, the constant $C$ may vary from line to line.

A permutation $\Gamma$ of real number set $S:=\{1, 2, \cdots, K_0\}$ is defined as a bijection  from $S$  to itself.  The cycle notation for permutation  expresses the permutation as a product of cycles corresponding to the orbits of the permutation.  So a permutation $\Gamma$ can be written as the product of disjoint cycles corresponding to the orbits   of the permutation, such as $\Gamma=\prod \mathcal{C}_l,$ where $\mathcal{C}_l$ is a cycle of length $l$.  Another fact is that any two pair partitions $\mathcal{P}_{2K_0}$ and     $\mathcal{P}_{2K_0}'$ form a permutation $\Gamma=\prod\mathcal{C}_l$ in $\{1,\cdots, 2 K_0\}$ and $\sum_{1\leq i_1\cdots , i_{2K_0}\leq p}\prod_{\{s,t\}\in\mathcal{P}_{2K_0}}\sigma_{i_si_t}\prod_{\{g,h\}\in\mathcal{P}'_{2K_0}}\sigma_{i_gi_h}=\prod_{\mathcal{C}_l\in \Gamma}
 \mathrm{Trace}(\Sigma^{|\mathcal{C}_l}).$ For example $\mathcal{P}_{4}=\{\{1,2\}, \{3,4\}\}$ and  $\mathcal{P}_{4}'=\{\{1,3\}, \{2,4\}\}$ form a permutation $\Gamma=(1,2,4,3)$ and  $\sum_{1\leq i_1\cdots , i_{4}\leq p}\prod_{\{s,t\}\in\mathcal{P}_{4}}\sigma_{i_si_t}\prod_{\{g,h\}\in\mathcal{P}'_{4}}\sigma_{i_gi_h}=
 \mathrm{Trace}(\Sigma^{|4}).$

\begin{IEEEproof}[Proof of Lemma \ref{iran}] By Lemma \ref{IS} (see Isserlis' theorem, Withers 1985),  we know that $\mathbb{E}\prod_{k=1}^{K_0}z_{i_{2k-1}}z_{i_{2k}}=\sum\prod \sigma_{i_si_t}$, where the notation $\sum$   means summing over all distinct ways of partitioning $\{1, \cdots, 2K_0 \}$ into pairs. This yields $(4K_0-1)! !$ terms in the sum.
 Then
  \begin{eqnarray*}
  &&\mathbb{E}\prod_{k=1}^{K_0}(z_{i_{2k-1}}z_{i_{2k}}-\sigma_{i_{2k-1}i_{2k}})\\
  &=&
  \mathbb{E}\prod_{k=1}^{K_0}z_{i_{2k-1}}z_{i_{2k}}-\sum_{l=1}^{K_0}\sigma_{i_{2l-1}i_{2l}}\mathbb{E}\prod_{k=1, k\neq l}^{K_0}z_{i_{2k-1}}z_{i_{2k}}\\
  &&+\sum_{1 \leq l<f\leq K_0}\sigma_{i_{2l-1}i_{2l}}\sigma_{i_{2f-1}i_{2f}}\prod_{k=1, k\neq l,f }^{K_0}z_{i_{2k-1}}z_{i_{2k}}-\cdots\\
  &&+(-1)^{K_0-1}\prod_{k=1}^{K_0} \sigma_{i_{2k-1}i_{2k}}. \end{eqnarray*}
Using the inclusion-exclusion principle, the above formula excludes out all the products which contains $\ \sigma_{i_{2k-1}i_{2k}}$ from $\sum$. 
So we get
            \begin{eqnarray*}\lbl{pee}
 &&\mathbb{E}\prod_{k=1}^{K_0}(z_{i_{2k-1}}z_{i_{2k}}-\sigma_{j_{2k-1}j_{2k}})\nonumber\\
 &=&\sum_{\mathcal{P}\in \mathcal{P}_{2K_0}}\prod_{\{s,t\}\in \mathcal{P}}\sigma_{i_si_t} ,
   \end{eqnarray*}
    where
  $\mathcal{P}_{2K_0}$ is the set contains    all the distinct ways  partitioning $\{1, \cdots,  2K_0 \}$ into pairs excluding the ways partitioning $ 2k-1$ and $2k $ into a pair for any $k\leq K_0$, and $\mathcal{P}$ is the set containing all the sub-index $\{s,t\}$ of  $\sigma_{i_si_t}$  which forms a way of pair partition in $\mathcal{P}_{2K_0}.$
Now let $j_{2k-1}=j_{2k}=i_{k}$ for all $1\leq k\leq K_0$, then
\begin{eqnarray}\lbl{poo}
&&\mathbb{E}\prod_{k=1}^{K_0}(z_{i_{k}}^2-\sigma_{i_ki_{k}})\nonumber \\
&=&\mathbb{E}\prod_{k=1}^{K_0}(z_{j_{2k-1}}z_{j_{2k}}-\sigma_{j_{2k-1}j_{2k}})\nonumber \\
&=& \sum_{\mathcal{P}\in \mathcal{P}_{2K_0}}\prod_{\{s,t\}\in \mathcal{P}}\sigma_{i_{\lceil \frac{s}{2}\rceil} i_{\lceil \frac{t}{2}\rceil}}.\end{eqnarray}
Moreover, we know that,  $\{s,t\} \in \mathcal{P}$, $\{s, t\}\neq \{2k-1,2k\}$ for any $k\leq K_0$, and $\cup_{ \mathcal{P}}\{s, t\}=\{1,\cdots, 2K_0\}$.
So  $\lceil \frac{s}{2}\rceil\neq\lceil \frac{t}{2}\rceil$, the all the sub-index $\{\lceil \frac{s}{2}\rceil ,\lceil \frac{t}{2}\rceil\}$ of $\sigma_{i_{\lceil \frac{s}{2}\rceil} i_{\lceil \frac{t}{2}\rceil}}$ in $\prod_{\mathcal{P}}$  after suitable position exchange, can form a bijection $\Gamma$ from $\{1,2, \cdots, K_0\}$ to itself.  So $\Gamma$ is a permutation of $\{1,\cdots, K_0\}$  without fixed points.  Then
    \begin{eqnarray*}\lbl{poo2}
 &&\sum_{1\leq i_1, \cdots, i_{K_0}\leq p}\mathbb{E}\prod_{k=1}^{K_0}(z_{1i_{k}}^2-\sigma_{i_ki_{k}})\nonumber \\
 &=&\sum_{\mathcal{P}\in \mathcal{P}_{2K_0}} \sum_{1\leq i_1, \cdots, i_{K_0}\leq p}\prod_{\{s,t\}\in \mathcal{P}}\sigma_{i_{\lceil \frac{s}{2}\rceil}i_{\lceil \frac{t}{2}\rceil}} \nonumber \\
 &\leq &Cp^{K_0/2}.
   \end{eqnarray*}
\end{IEEEproof}

\begin{IEEEproof}[Proof of Lemma \ref{rain}] The first equality in (\ref{tian}) is due to the definition of the trace of a matrix. By the assumption of $\lambda_{max}(\Sigma)\leq M_0$, we induce that $\sigma_{ii}\leq M_0$ for $i \geq 1$, and $|\sigma_{ij}|$ have a universal bound. Then by Perron--Frobenius Theorem and $\Sigma\in \mathcal{F}_{\alpha}$, we can set that $\lambda_{max}(\Sigma_a)\leq \max_{j}\sum_{i}|\sigma_{ij}|\leq 3M_0+M$. Then   $\mathrm{Trace}(\Sigma_a^k)\leq p(3M_0+M) ^k$.
\end{IEEEproof}

\begin{IEEEproof}[Proof of Lemma \ref{day}]  Recall the definitions of $a_n. b_n, a_{ij}, b_{ij}$ and $ \omega_{ij}{(\tau)}$,  then
 \begin{eqnarray*}
  \bar{A}^{(\tau)}_ {ij}&=&\frac{n^2}{(n-1)^2}\bigg\{\frac{(C_n-1)(n-3)}{(n-2)(n+1)}-\frac{C_n^2(n-3)^2}{(n-2)^2(n+1)^2}\bigg\}\nonumber \\&&+\bigg[\omega_{ij} ^{(\tau)}-\frac{n}{n-1}\bigg\{1-\frac{C_n(n-3)}{2(n-2)(n+1)}\bigg\}\bigg]^2,\nonumber\\ \, \,\bar{B}^ {(\tau)}_ {ij}&=&  (\omega_{ij} ^{(\tau)})^2+\frac{n}{n-1}(C_n-2) \omega_{ij} ^{(\tau)}, \nonumber\\          \bar{C}^{(\tau)}_ {ij}&=&\frac{n^2}{(n-1)^2}\bigg\{\frac{2(C_n-1)}{n+1}-\frac{C_n^2}{(n+1)^2}\bigg\}\nonumber \\
  &&+\bigg\{\omega_{ij} ^{(\tau)}-\frac{n }{n-1}(1-\frac{C_n}{n+1})\bigg\}^2  \lbl{fan}.
      \end {eqnarray*}
       So  $\bar{A}^ {(\tau)}_{ij}\leq \frac{n^2}{(n-2)(n+1)}$, $\bar{B}^{(\tau)}_ {ij}\leq C_n$ and $\bar{C}^{(\tau)}_ {ij}\leq\frac{n^2}{n^2-1} $.
     It is easy to see that $\bar{A}^ {(\tau)}_{ij}$ and $\bar{C}^{(\tau)}_ {ij}$ are bounded by 2. Furthermore, $ \bar{B}^ {(\tau)}_ {ij}=0$ when $|i-j|\geq \tau$ and $\bar{B}^{(\tau)}_ {ij}
\geq \frac{C_n}{2}$ for any $|i-j|\leq \lfloor\frac{\tau}{2}\rfloor$. When n is large enough, there exists $\epsilon>0$ such that $ \bar{A}^ {(\tau)}_ {ij}\leq1+\epsilon$, $\bar{C}^{(\tau)}_ {ij}\leq 1+\epsilon$.
\end{IEEEproof}

 In the following proofs, we  define that   $\bar{A}^ {(\tau)}_ {i_1\cdots i_{2k}}=\prod_{j=1}^k\bar{A}^ {(\tau)}_ {i_{2j-1}i_{2j}}$   and   $\bar{B}^ {(\tau)}_ {i_1\cdots i_{2k}}=\prod_{j=1}^k\bar{B}^ {(\tau)}_ {i_{2j-1}i_{2j}}$.  Let  $\mathbb{E}_{k-1}$ denotes $\mathbb{E}(\cdot |\ml{F}_{k-1})$.

\begin{IEEEproof}[Proof of Lemma \ref{zhu}]  First by the definitions of $Y_{nm}$, $Y_m$ and $U_m$,
we observe that
\begin{eqnarray*}\lbl{zui} && \mathbb{E}_{m-1}Y_{nm}^2\nonumber \\
   &=&\frac{1}{n^4}\sum_{l,f=1}^{m-1} \sum_{1\leq i, j,s,t\leq p}[\nonumber \\
   &&4\bar{A}^ {(\tau)}_ {ij}\bar{A}^ {(\tau)}_ {st}(\sigma_{is}\sigma_{jt}+\sigma_{it}\sigma_{js})  (z_{li}z_{lj}-\sigma_{ij})(z_{fs}z_{ft}-\sigma_{st})]\nonumber\\
   && \mathbb{E}_{m-1}Y_{nm}Y_m \nonumber \\
      &=&\frac{1}{n^4}\sum_{l=1}^{m-1}\sum_{1\leq i, j,s,t  \leq p} [ \nonumber \\
     && 4 (n-2)\bar{A}^ {(\tau)}_ {ij} \bar{A}^ {(\tau)}_ {st}\sigma_{st}(\sigma_{is}\sigma_{jt}+\sigma_{it}\sigma_{js})(z_{li}z_{lj}-\sigma_{ij})]\nonumber\\
    &&\mathbb{E}_{m-1}Y_{nm}U_m \nonumber \\
   &=&\frac{1}{n^4}\sum_{l=1}^{m-1}\sum_{1\leq i, j,s,t  \leq p}[ \nonumber \\
  && 4\bar{A}^ {(\tau)}_ {ij} \bar{B}^ {(\tau)}_ {st}(\sigma_{ss}\sigma_{it} \sigma_{jt}+\sigma_{tt}\sigma_{is}\sigma_{js})(z_{li}z_{lj}-\sigma_{ij})].\end{eqnarray*}
 Applying Lemmas \ref{iran}, \ref{rain} and \ref{day}, we can calculate that
 \begin{eqnarray*}\lbl{group}
    &&\sum_{m=2}^{n-1}\sum_{k=2}^{n-1}\mathbb{E}[\mathbb{E}_{m-1}(Y_{nm}Y_m)\mathbb{E}_{k-1}(Y_{nk}Y_k)]\nonumber\\
    &\leq& \frac{C}{n^3}\mathrm{Trace}(
   \Sigma_a^8) =C\frac{p}{n^3}         \end{eqnarray*}
  and
   \begin{eqnarray*}\lbl{diss}
   &&\sum_{m=2}^{n-1}\sum_{k=2}^{n-1}\mathbb{E}[\mathbb{E}_{m-1}Y_{nm}^2\mathbb{E}_{k-1}(Y_{nk}U_k)]\nonumber\\
    &\leq& \frac{CC_n}{n^5}2\tau M_0\mathrm{Trace}(\Sigma^2)\mathrm{Trace}(\Sigma_a^5)=C\frac{p^2\tau C_n}{n^5}.
  \end{eqnarray*}
Then by the similar arguments as above, as $n\to \infty$, we have that
   \begin{eqnarray*}
   &&\sum_{m=2}^{n-1}\sum_{k=2}^{n-1}\mathbb{E}[\mathbb{E}_{m-1}Y_{nm}^2\mathbb{E}_{k-1}(Y_{nk}Y_k)] \\
& \leq & \frac{C}{n^4}\mathrm{Trace}(\Sigma^2)\mathrm{Trace}(\Sigma_a^6)=C\frac{p^2}{n^4};\nonumber\\
  &&\sum_{m=2}^{n-1}\sum_{k=2}^{n-1}\mathbb{E}[\mathbb{E}_{m-1}(Y_{nm}Y_m)\mathbb{E}_{k-1}(Y_{nk}U_k)]\\
  &   \leq & \frac{CC_n}{n^4}2\tau M_0\mathrm{Trace}(
   \Sigma_a^7) =C\frac{p\tau C_n}{n^4};\nonumber\\
    &&\sum_{m=2}^{n-1}\sum_{k=2}^{n-1}\mathbb{E}[\mathbb{E}_{m-1}(Y_{nm}U_m)\mathbb{E}_{k-1}(Y_{nk}U_k)] \\
    & \leq & \frac{CC_n^2}{n^5}(2\tau M_0)^2\mathrm{Trace}(
   \Sigma_a^6) =C\frac{p(\tau C_n)^2}{n^5}.        \end{eqnarray*}
    So we can conclude that $D_n\leq C(\frac{p(\tau C_n)^2}{n^5}+\frac{p^2 \tau C_n}{n^5}+\frac{p^2}{n^4})$.
    \end{IEEEproof}

\begin{IEEEproof}[Proof of Lemma \ref{son}]
     Following the same notation as in the proof of Theorem \ref{girl}, we have that
       \begin{eqnarray*}s_{n-1}^2 =\sum_{m=2}^{n-1}\mathbb{E}[\mathbb{E}_{m-1}Y_{nm}^2]+(n-1) (\mathbb{E} Y_{m}^2+\mathbb{E} U_{m}^2+2\mathbb{E} Y_{m}U_m)\nonumber
              \end{eqnarray*}
where          \begin{eqnarray*}
 &&\mathbb{E}Y_{nm}^2\\
  &=&\frac{1}{n^4}(m-1)\sum_{1\leq i, j,s,t\leq p}4\bar{A}^ {(\tau)}_ {ij}\bar{A}^ {(\tau)}_ {st}(\sigma_{is}\sigma_{jt}+\sigma_{it}\sigma_{js}) ^2\\
  &\leq &\frac{Cmp^2}{n^4} ; \\
  && \mathbb{E}Y_m^2\\
  &=&  \frac{4(n-2)^2}{n^4} \sum_{1\leq i, j ,s,t \leq p}\bar{A}^ {(\tau)}_ {ij}\bar{A}^ {(\tau)}_ {st}\sigma_{ij}\sigma_{st}(\sigma_{is}\sigma_{jt}+\sigma_{it}\sigma_{js})\\
  &\leq& C\frac{p}{n^2};
   \end{eqnarray*}
 \begin{eqnarray*}
   &&\mathbb{E}Y_mU_m \nonumber  \\
   &=&  \frac{4(n-2) }{n^4} \sum_{1\leq i, j ,s,t \leq p}\bar{A}^ {(\tau)}_ {ij}\bar{B}^ {(\tau)}_ {st}\sigma_{ij}( \sigma_{ss}\sigma_{it} \sigma_{jt} \nonumber \\
   &&+\sigma_{tt}\sigma_{is}\sigma_{js});\lbl{tea}\\ &&\mathbb{E}U_m^2 \nonumber \\
      &=&\frac{1}{n^4} \sum_{1\leq i, j,s,t \leq p}2\bar{B}^ {(\tau)}_ {ij}\bar{B}^ {(\tau)}_ {st} \\
      &&(\sigma_{ii}\sigma_{ss}\sigma_{jt}^2+\sigma_{ii}\sigma_{tt}\sigma_{js}^2+\sigma_{jj}\sigma_{ss}\sigma_{it}^2+\sigma_{jj}\sigma_{tt}\sigma_{is}^2).\nonumber  \end{eqnarray*}
      So it is easy to check that $s_{n-1}^2$  has the form as in (\ref{var}).
Since $\Sigma$ satisfies the condition (C.2),  $\sigma_{ij}\leq M_0 $ for any $1\leq i, j\leq p$, and with (\ref{tian}) we have
  \begin{eqnarray}
&&\sum_{  |i-j|< \tau,|s-t|< \tau}  \sigma_{ii}\sigma_{ss}\sigma_{jt}^2 \nonumber\\
&=&\sum_{1\leq j, t \leq p}\sigma_{jt}^2 (\sum_{i: |i-j|< \tau}\sigma_{ii} ) (\sum_{s: |s-t|< \tau} \sigma_{ss})\nonumber\\
&\leq& (2\tau M_0)^2\mathrm{Trace}(\Sigma^2)\leq 4M_0^3p\tau^2.\lbl{kiss}\\
 &&\sum_{ 1 \leq i, j\leq p ,|s-t|< \tau}  \sigma_{ij}\sigma_{ss}\sigma_{it} \sigma_{jt}\nonumber \\
  &=& \sum_{t=1}^p (\sum_{s: |s-t|< \tau}\sigma_{ss} ) [\sum_{i=1}^p\sum_{j=1}^p\sigma_{ij}\sigma_{jt}\sigma_{it}] \nonumber\\
  & \leq& 2M_0p \tau.\lbl{kick}
\end{eqnarray}
From Lemmas \ref{day} and {\ref{rain}},  we know that $\bar{B}^ {(\tau)}_ {ij}=0$ when $|i-j| \geq \tau$, and  (\ref{tian}), it is easy to see that
  \begin{eqnarray*}\lbl{like}  &&\mathbb{E}U_m^2 \nonumber \\
  &\leq&   \frac{2C_n^2}{n^4}\sum_{  |i-j|< \tau,|s-t|< \tau}  (\sigma_{ii}\sigma_{ss}\sigma_{jt}^2+\sigma_{ii}\sigma_{tt}\sigma_{js}^2+\sigma_{jj}\sigma_{ss}\sigma_{it}^2 \nonumber \\
  &&+\sigma_{jj}\sigma_{tt}\sigma_{is}^2) \nonumber \\
  &\leq&  C \frac{p \tau^2C_n^2}{n^4} \quad \end{eqnarray*}
  and let $i$ or $j$ equal $t$ or $s$ in the above summation
  \begin{eqnarray*}\lbl{like1}  \mathbb{E}U_m^2 &\geq&  \frac{2C_n^2}{n^4}\sum_{  |i-j|\leq \lfloor\frac{\tau}{2}\rfloor, |s-j|\leq \lfloor\frac{\tau}{2}\rfloor}  \sigma_{ii}\sigma_{ss}\sigma_{jj}^2.\end{eqnarray*}
Since there exists $C$ such that $\min_{1\leq i\leq p}\sigma_{ii}>\frac{2}{C^{1/4}}$, $ \mathbb{E}U_m^2\geq \frac{p\tau^2C_n^2}{Cn^4} $.
 Furthermore, using (\ref{kick}), we have
  \begin{eqnarray*}\lbl{liked}  &&\mathbb{E}Y_mU_m \nonumber \\
  &\leq &  \frac{4(n-2)C_n}{n^4}\sum_{ 1\leq i,j\leq p,  |s-t|<  \tau }  \nonumber \\
  && \sigma_{ij}( \sigma_{ss}\sigma_{it} \sigma_{jt}+\sigma_{tt}\sigma_{is}\sigma_{js})\\
 &\leq & C\frac{p\tau C_n}{n^3}  .\end{eqnarray*}
 By Lemma \ref{day}, we also know that   $|\bar{A}^ {(\tau)}_ {ij}|$ is bounded by 2 for $1\leq i, j\leq p$, then
                 \begin{eqnarray*}s_{n-1}^2
  \leq C\big(\frac{p^2}{n^2}+\frac{  p\tau^2C_n^2}{n^3} +\frac{p\tau C_n}{n^2}  \big).\lbl{up}
                 \end{eqnarray*}
Since $|\bar{A}^ {(\tau)}_ {ij}|\geq 1$ for all $|i-j|>\tau $, we have
 \begin{eqnarray*}
 &&\sum_{m=2}^{n-1}\mathbb{E}[\mathbb{E}_{m-1}Y_{nm}^2]\\
 &\geq & \frac{2(n-1)(n-2)}{n^4} \sum_{| i-j| > \tau}(\bar{A}^ {(\tau)}_ {ij})^2 (\sigma_{ii}\sigma_{jj}+\sigma_{ij}^2)^2 \\
 & \geq & \frac{(p-\tau)^2}{Cn^2})
   \end{eqnarray*}  and $\mathbb{E}U_m^2\gg \mathbb{E}Y_m^2$ when $\tau\gg \frac{n}{C_n}$. Hence,
\begin{eqnarray}\lbl{chi}s_{n-1}^2&\geq&\sum_{m=2}^{n-1}\mathbb{E}[\mathbb{E}_{m-1}Y_{nm}^2]+(n-1) (\sqrt{\mathbb{E} Y_{m}^2}-\sqrt{\mathbb{E} U_{m}^2})^2 \nonumber\\
&\geq&
\left\{
  \begin{array}{l l } \frac{1}{C}(\frac{(p-\tau)^2}{n^2}+\frac{p\tau^2C_n^2}{n^3})
    & \quad \text{when $ \tau\gg \frac{n}{C_n}$};\\
   \frac{1}{C}\frac{(p-\tau)^2}{n^2}   & \quad \text{otherwise.}
  \end{array} \right.\
  \end{eqnarray}
 So we can conclude that $s_{n-1}^2\geq \frac{1}{C} \max (\frac{p^2}{n^2},\frac{p\tau^2C_n^2}{n^3} )$ when $C_n$ is not a constant and $C_n\to \infty $ as $n\to \infty$ or $p\gg n$.
  When $C_n$ is a constant and $\frac{p}{n}\to y < \infty$, then $s_{n-1}^2\geq \frac{1}{C}\frac{p^2}{n^2} $ for any $\epsilon>0$ and $\tau\leq (1-\epsilon)p$.
  \end{IEEEproof}

\begin{IEEEproof}[Proof of Lemma \ref{night}] By the definition of $H_n(Z_1, Z_2)$ and Lemma \ref{iran}, it is easy to see that, for any even number $K_0$,
\begin{eqnarray*} \lbl{sweet}
  && \mathbb{E}H_n(Z_1, Z_2)^{K_0}\\
  &=&\frac{2^{K_0}}{n^{2K_0}} \sum_{ 1\leq i_1, \cdots, i_{2K_0} \leq p} \bar{A}^ {(\tau)}_ {i_1\cdots i_{2K_0}}   [\\
  &&\mathbb{E}\prod_{k=1}^{K_0}(z_{1i_{2k-1}}z_{1i_{2k}}-\sigma_{i_{2k-1}i_{2k}})   ]^2  \nonumber\\ &\leq& \frac{C2^{K_0}}{n^{2K_0}}\sum_{\mathcal{P}\in \mathcal{P}_{2K_0}} \sum_{\mathcal{P}'\in \mathcal{P}_{2K_0}} ( \\
  &&\sum_{ 1\leq i_1, \cdots, i_{2K_0} \leq p} \prod_{\{s,t\}\in \mathcal{P}}\prod_{\{g,h\}\in \mathcal{P}'}\sigma_{i_si_t}\sigma_{i_gi_h}),\end{eqnarray*}  where $\mathcal{P}_{2K_0}$ is defined as in Lemma \ref{iran} and $\mathcal{P}$, $\mathcal{P}'\in \mathcal{P}_{2K_0}$.  Note that all the  sub-index $\{u, v\}$ of $ i_ui_v$ from $\sigma_{i_ui_v}$ in the product  $\prod_{\{s,t\}\in \mathcal{P}}\prod_{\{g,h\}\in \mathcal{P}'}\sigma_{i_si_t}\sigma_{i_gi_h}$, after suitable position exchange, can form a bijection $\Gamma$ from $\{1,2, \cdots, 2K_0\}$ to itself.
So $\Gamma$ is a permutation of $\{1,\cdots, 2K_0\}$ without 1-element cycle. Then using cycle notation, $\Gamma=\prod_{l}\mathcal{C}_l$, which $\mathcal{C}_l$ are the disjointed cycle with length $|\mathcal{C}_l|=l\geq 2$ forms the permutation $\Gamma$ and $\sum_{l}|\mathcal{C}_l|=2K_0$, so
\begin{eqnarray*}
&&\sum_{ 1\leq i_1, \cdots, i_{2K_0} \leq p} \prod_{\{s,t\}\in \mathcal{P}}\prod_{\{g,h\}\in \mathcal{P}'}|\sigma_{i_si_t}\sigma_{i_gi_h}|\\
&=& \prod_{\mathcal{C}_l}\mathrm{Trace}(\Sigma_a^l)\\
&\leq& Cp^{K_0}.
\end{eqnarray*} So  $ \mathbb{E}H_n(Z_1, Z_2)^{K_0} \leq  C\frac{p^{K_0}}{n^{2K_0}}$.

Now we calculate the $K_0$ th moment of  $Y_m$,
  \begin{eqnarray*} \lbl{sweet}
    &&\mathbb{E}Y_{m}^{K_0}\\
  &=& \frac{2^{K_0} (n-2)^{K_0}}{n^{2K_0}}  \sum_{1\leq i_1, \cdots, i_{2K_0} \leq p} \nonumber \\
  &&\bar{A}^ {(\tau)}_ {i_1\cdots i_{2K_0}} \sigma_{i_1i_2}\cdots \sigma_{i_{2K_0-1}i_{2K_0}}  ( \sum_{\mathcal{P}\in \mathcal{P}_{2K_0}}\prod_{\{s,t\}\in \mathcal{P}}\sigma_{i_si_t}) \nonumber\\
&\leq&C  \frac{2^{K_0}}{n^{K_0}} \sum_{\mathcal{P}\in \mathcal{P}_{2K_0}}\\
&&( \sum_{1\leq i_1, \cdots, i_{2K_0} \leq p} |\sigma_{i_1i_2}\cdots \sigma_{i_{2K_0-1}i_{2K_0}}|\prod_{\{s,t\}\in \mathcal{P}}|\sigma_{i_si_t}|).\end{eqnarray*} Using the similar arguments as  above,  the sub-index  $\{s, t\}$ of $i_si_t$ from $\sigma_{i_si_t}$ in the product $\sigma_{i_1i_2}\cdots \sigma_{i_{K_0-1}i_{K_0}}|\prod_{(s,t)\in \mathcal{P}}|\sigma_{i_si_t}|$ generates 
a  permutation $\Gamma'=\prod_l \mathcal{C}_l$ by using cycle notation. Furthermore all $\mathcal{C}_l$ are the cycle with length $|\mathcal{C}_l|\geq 4$. This is because that   $\{s,t\}\in \mathcal{P}  \neq \{2k-1, 2k\}$, then there  is no cycle with length 2 in the permutation $\Gamma'$. If  there is a cycle $\mathcal{C}_3$ with length  3, without loss of generality,   assume this cycle starts with $k$ (an odd number), so $\mathcal{C}_3$ can be written as $k\mapsto k+1\mapsto l\mapsto k$, while $l\neq k\neq k+1$. Since $\{k, k+1\}$ is the sub-index from the terms $\sigma_{i_1i_2}\cdots \sigma_{i_{2K_0-1}i_{K_0}}$ and $l\neq k$, then $\{k+1,l\}$ comes from $\mathcal{P}$ and $\{l,k\}\in\{\{2m-1,2m\}, 1\leq m \leq K_0\}$, which is a contradiction to $l\neq k+1.$ So the smallest length of $\mathcal{C}_l$ is not smaller than 4. Then $ \mathbb{E}Y_{m}^{K_0} \leq C\frac{p^\frac{K_0}{2}}{n^{K_0}}.$

Using the similar argument  for (\ref{poo}) as in Lemma \ref{iran}, we get that
\begin{eqnarray*}
 && \mathbb{E}U_{m}^{K_0}\\
 &=& \frac{2^{K_0}}{n^{2K_0}}   \sum_{\substack{ |i_{2j-1}-i_{2j}|\leq \tau \\ 1 \leq j\leq K_0}} \\
  &&\bar{B}^ {(\tau)}_ {i_1\cdots i_{2K_0}}\sigma_{i_1i_1}\sigma_{i_3i_3}\cdots \sigma_{i_{2K_0-1}i_{2K_0-1}}  \mathbb{E}\prod_{k=1}^{K_0}(z_{1i_{2k}}^2-\sigma_{i_{2k}i_{2k}})\nonumber\\
   &\leq&  \frac{(2C_n)^{K_0}}{n^{2K_0}}  \sum_{\mathcal{P}\in \mathcal{P}_{2K_0}}  \sum_{\substack{ |i_{2j-1}-i_{2j}|\leq \tau \\ 1\leq j\leq K_0}}  \\
    && \sigma_{i_1i_1}  \sigma_{i_3i_3}\cdots \sigma_{i_{2K_0-1}i_{2K_0-1}}  \prod_{\{s,t\}\in \mathcal{P}}|\sigma_{i_{2\lceil \frac{s}{2}\rceil} i_{2\lceil \frac{t}{2}\rceil}}|.\nonumber
    \end{eqnarray*}
    Then with the property of $\mathcal{P}$, we know that $2\lceil \frac{s}{2}\rceil\neq2\lceil \frac{t}{2}\rceil$. All the sub-index $\{2\lceil \frac{s}{2}\rceil ,2\lceil \frac{t}{2}\rceil\}$ of $\sigma_{i_{2\lceil \frac{s}{2}\rceil} i_{2\lceil \frac{t}{2}\rceil}}$ in $\prod_{\mathcal{P}}$ form a permutation $\Gamma_{\mathcal{P}}$ of $\{2,4, \cdots, K_0\}$ without fixed points, while $\Gamma_{\mathcal{P}}=\prod \mathcal{C}_l$ and $\mathcal{C}_l$ is a cycle with length $l\geq 2$, so there are at most $\frac{K_0}{2}$ cycles in the permutation $\Gamma$. Further $\mathcal{C}_l=(j_1, j_2, \cdots, j_l) $, where $j_i\in \{2, 4, \cdots, K_0\}$. Then
    \begin{eqnarray*}
  &&\mathbb{E}U_{m}^{K_0} \\
   &\leq&  \frac{(2C_n)^{K_0}}{n^{2K_0}}  \sum_{\mathcal{P}\in \mathcal{P}_{2K_0}}  \prod_{\mathcal{C}_l\in \Gamma_{\mathcal{P}}}\\
   &&\{\sum_{\substack{ |i_{(j_h-1)}-i_{j_h}|< \tau \\ 1\leq h\leq l}}  (\prod_{k=1}^l\sigma_{i_{(jk-1)}i_{(jk-1)}} )  |\sigma_{i_{j_1} i_{j_2}}\cdots \sigma_{i_{j_l}i_{j_1}}|\}\nonumber\\
   &\leq&  \frac{(2C_n)^{K_0}}{n^{2K_0}}  \sum_{\mathcal{P}\in \mathcal{P}_{2K_0}}  \prod_{\mathcal{C}_l\in \Gamma_{\mathcal{P}}} \\
   &&[2\tau M_0]^l(\sum_{\substack{ 1 \leq i_{j_h}\leq p\\ 1\leq h\leq l}}   |\sigma_{i_{j_1} i_{j_2}}\cdots \sigma_{i_{j_l}i_{j_1}}|)\nonumber\\
  &\leq&  \frac{(4 M_0 C_n\tau)^{K_0}}{n^{2K_0}} \sum_{\mathcal{P}\in \mathcal{P}_{2K_0}}  \prod_{\mathcal{C}_l\in \Gamma_{\mathcal{P}}}Trace(\Sigma_a^l)\\
  &\leq& C \frac{(C_n\tau)^{K_0}p^{\frac{K_0}{2}}}{n^{2K_0}}.
         \end{eqnarray*} Now we get all the results.
\end{IEEEproof}

\begin{IEEEproof}[Proof of Lemma \ref{dress}]  Define
\begin{eqnarray*}
&&T_n(Z_l, Z_f)\\
&=&\mathbb{E}_{m-1} H_n(Z_m, Z_l)H_n(Z_m, Z_f)  \nonumber\\
&=&\frac{1}{n^4} \sum_{1\leq i, j,s,t\leq p}4\bar{A}^ {(\tau)}_ {ij}\bar{A}^ {(\tau)}_ {st}\\
&&(\sigma_{is}\sigma_{jt}+\sigma_{it}\sigma_{js})(z_{li}z_{lj}-\sigma_{ij})(z_{fs}z_{ft}-\sigma_{st})\lbl{kim},
\end{eqnarray*}  where $l, f\leq m-1$ for all $m \geq 2$. Then  \begin{eqnarray*} &&\mathbb{E}_{m-1}Y_{nm}^2\\
&=&\sum_{l=1}^{m-1}\sum_{f=1}^{m-1} T_n(Z_l, Z_f)\\
&=&\sum_{l=1}^{m-1}T_n(Z_l, Z_l)+2\sum_{1\leq l < f \leq m-1 }T_n(Z_l, Z_f)\end{eqnarray*}
and $\mathbb{E}Y_{nm}^2=\sum_{l=1}^{m-1}\mathbb{E}T_n(Z_l, Z_l)$. So by  the moments' bounds for Martingale and $C_r$ inequality,
         \begin{eqnarray}\lbl{destroy}
 &&\mathbb{E}[\sum_{m=2}^{n-1}(\mathbb{E}_{m-1}Y^2_{nm}-\mathbb{E}Y_{nm}^2)]^{K_0} \nonumber \\
  &\leq& 2^{K_0-1}\{L_{K_0}n^{\frac{3K_0}{2}}\mathbb{E}[ T_n(Z_1,Z_1)-\mathbb{E}T_n(Z_1, Z_1)]^{K_0}\nonumber\\
  &&+L_{K_0}^2n^{2K_0}\mathbb{E}T_n(Z_1, Z_2)^{K_0}\}.
          \end{eqnarray}
    First we show that $\mathbb{E}[ T_n(Z_1,Z_1)-\mathbb{E}T_n(Z_1, Z_1)]^{K_0}\leq C\frac{p^{3K_0/2}}{n^{4K_0}}$. Let $E_l=(\sigma_{i_{(4l-3)}i_{(4l-1)}}\\\sigma_{i_{(4l-2)}i_{4l}}+\sigma_{i_{(4l-3)}i_{4l}}\sigma_{i_{(4l-1)}i_{(4l-2)}})$, we have
    \begin{eqnarray*} && \mathbb{E}[ T_n(Z_1,Z_1)-\mathbb{E}T_n(Z_1, Z_1)]^{K_0}\nonumber\\&&= \frac{4^{K_0}}{n^{4K_0}} \sum_{ 1\leq i_1, \cdots, i_{4K_0}\leq p} \bar{A}^ {(\tau)}_ {i_1\cdots i_{4K_0}} \prod_{l=1}^{K_0}\{ \\
    &&E_l\times[\mathbb{E}(z_{1i_{(4l-3)}}z_{1i_{(4l-2)}}-\sigma_{i_{(4l-3)}i_{(4l-2)}})\nonumber\\
    &&(z_{1i_{(4l-1)}}z_{1i_{4l}}-\sigma_{i_{(4l-1)}i_{4l}})-E_l]\}.\end{eqnarray*}
      Using the same inclusion and exclusion principle   argument as in Lemma \ref{iran}, we conclude that   \begin{eqnarray*} &&\prod_{l=1}^{K_0}\mathbb{E}[(z_{1i_{(4l-3)}}z_{1i_{(4l-2)}}-\sigma_{i_{(4l-3)}i_{(4l-2)}})\\
      &&(z_{1i_{(4l-1)}}z_{1i_{4l}}-\sigma_{i_{(4l-1)}i_{4l}})-E_l]\nonumber
      \\&=&\sum_{\mathcal{P}\in\mathcal{P}_{4K_0}}\prod_{\{s,t\}\in \mathcal{P}}\sigma_{st},\end{eqnarray*} where $\mathcal{P}_{4K_0}$ is   the set containing   all the distinct ways  that partition $1, \cdots,  4K_0 $ into pairs excluding the ways that partition $ 2k-1$ and $2k $ into a pair for any $k\leq 2K_0$ and the ways that contain $\{4l-3, 4l-1\}$  and $\{4l-2,4l\}$ or $\{4l-3 ,4l\}$ and $\{4l-2,4l-3\}$  for any $l\leq K_0$ and $\mathcal{P} \in \mathcal{P}_{4K_0}$. By the similar argument as in the proof of Lemma \ref{night}, all the sub-index $\{g,h\}$ from $\sigma_{i_gi_h}$ in the product $[\prod_{l=1}^{K_0} (|\sigma_{i_{(4l-3)}i_{(4l-1)}}||\sigma_{i_{(4l-2)}i_{4l}}|+|\sigma_{i_{(4l-3)}i_{4l}}||\sigma_{i_{(4l-1)}i_{(4l-2)}}|)] \times\prod_{\{s,t\}\in \mathcal{P}}|\sigma_{i_si_t}| $ form a permutation $\Gamma$ in $\{1, \cdots, 4K_0\}$.   By the definition of $\mathcal{P}$, we know there are at most $K_0$ cycles with length 2 in $\Gamma$.  We claim there is no cycle with length 3 in $\Gamma$. If there exits such a cycle with length 3, without loss of generality, we assume it can be written as $4l-3\mapsto 4l-1\mapsto s\mapsto 4l-3$.  Whence $\{4l-1,s\}\in \mathcal{P}$ and $\{ 4l-3, s\}\in \{\{2t-1, 2t\}, t\leq 2K_0\}$ then $s=4l-1$, then there is a contradiction. So   $\frac{3K_0}{2}$ is the largest number of cycles to form  such a  permutation $\Gamma$. In this case, there are  $K_0$ cycles with length 2 and $\frac{K_0}{2}$ cycles with length 4 that give the largest upper bound for
  \begin{eqnarray*}&&\sum_{ 1\leq i_1, \cdots, i_{4K_0}\leq p}[\prod_{l=1}^{K_0} (|\sigma_{i_{(4l-3)}i_{(4l-1)}}||\sigma_{i_{(4l-2)}i_{4l}}|\\
  &&+|\sigma_{i_{(4l-3)}i_{4l}}||\sigma_{i_{(4l-2)}i_{(4l-3)}}|)] \times\prod_{\{s,t\}\in \mathcal{P}}|\sigma_{st}|\nonumber
  \\&&\leq \mathrm{Trace}(\Sigma_a^2)^{K_0}\mathrm{Trace}(\Sigma_{a}^4)^{\frac{K_0}{2}}= Cp^{3K_0/2}.
  \end{eqnarray*} So we get that
   \begin{eqnarray}\lbl{dinner}\mathbb{E}[ T_n(Z_1,Z_1)-\mathbb{E}T_n(Z_1, Z_1)]^{K_0}\leq C\frac{p^{3K_0/2}}{n^{4K_0}}.\end{eqnarray}
   Now we need to find the bound for $\mathbb{E}T_n(Z_1, Z_2)^{K_0}$.   First
   \begin{eqnarray*} &&\mathbb{E}T_n(Z_1,Z_2) ^{K_0}\\
   &=& \frac{4^{K_0}}{n^{4K_0}} \sum_{ 1\leq i_1, \cdots, i_{4K_0}\leq p} \bar{A}^ {(\tau)}_ {i_1\cdots i_{4K_0}} \prod_{l=1}^{K_0}E_l\times \\
   &&\mathbb{E}(z_{1i_{(4l-3)}}z_{1i_{(4l-2)}}-\sigma_{i_{(4l-3)}i_{(4l-2)}})\nonumber\\
   &&\mathbb{E}(z_{1i_{(4l-1)}}z_{1i_{4l}}-\sigma_{i_{(4l-1)}i_{4l}}).\end{eqnarray*}
    Using  Lemma \ref{iran}, we  have
    $\prod_{l=1}^{K_0}\mathbb{E}(z_{1i_{(4l-3)}}z_{1i_{(4l-2)}}-\sigma_{i_{(4l-3)}i_{(4l-2)}})
     =\sum_{\mathcal{P}_{4K_0}}\prod_{\{s,t\}\in \mathcal{P}}\sigma_{i_si_t},
    $
    where $\mathcal{P}_{4K_0}$ is   the set contains   all the distinct ways  partitioning $1, 2, 5, 6,  \cdots,4K_0-3,  4K_0-2 $ into pairs excluding the ways partitioning $ 4l-3$ and $4l-2 $ into a pair for any $l\leq K_0$  and  $\mathcal{P} \in \mathcal{P}_{4K_0}$. Similarly,
$\prod_{l=1}^{K_0}\mathbb{E}(z_{1i_{(4l-1)}}z_{1i_{(4l)}}-\sigma_{i_{(4l-1)}i_{4l}}) =\sum_{\mathcal{P}'_{4K_0}}\prod_{\{s,t\}\in \mathcal{P}'}\sigma_{i_si_t},$
where $\mathcal{P}'_{4K_0}$ is   the set contains   all the distinct ways  that partitioning $3, 4, 7, 8,  \cdots,4K_0-1,  4K_0 $ into pairs excluding the ways  partitioning $ 4l-1$ and $4l $ into a pair for any $l\leq K_0$  and $\mathcal{P}'$ $\in \mathcal{P}'_{4K_0}$.
   By the similar argument as in the proof of Lemma \ref{night}, all the sub-index $\{u,v\}$ from $\sigma_{i_si_t}$ in the product $[\prod_{l=1}^{K_0} (|\sigma_{i_{(4l-3)}i_{(4l-1)}}||\sigma_{i_{(4l-2)}i_{4l}}|+|\sigma_{i_{(4l-3)}i_{4l}}||\sigma_{i_{(4l-1)}i_{(4l-2)}}|)] \prod_{\{s,t\}\in
   \mathcal{P}}|\sigma_{i_si_t}| \prod_{\{g,h\}\in
   \mathcal{P}'}|\sigma_{i_gi_h}| $ form a permutation $\Gamma'$ in $\{1, \cdots, 4K_0\}$. We claim there is no cycle with length 2 or 3 in $\Gamma'$. It is easy to check that there is no cycle in $\Gamma'$ with length 2.  If there exits a cycle with length 3, without loss of generality, we assume it can be written as $4l-3\mapsto 4l-1\mapsto s\mapsto 4l-3$.  Whence $\{4l-1,s\}\in \mathcal{P}'$ and $\{s, 4l-3\}\in \mathcal{P}$. However $\mathcal{P}_{4K_0}$ and $\mathcal{P}'_{4K_0}$ do not contain the same elements, then there is a contradiction. So we get that
  \begin{eqnarray*}&&\sum_{ 1\leq i_1, \cdots, i_{4K_0}\leq p}\prod_{l=1}^{K_0} (|\sigma_{i_{(4l-3)}i_{(4l-1)}}||\sigma_{i_{(4l-2)}i_{4l}}|\\
  &&+|\sigma_{i_{(4l-3)}i_{4l}}||\sigma_{i_{(4l-2)}i_{(4l-3)}}|) \prod_{\substack{\{s,t\}\in \mathcal{P}\\\{g,h\}\in \mathcal{P}'}}|\sigma_{st}| |\sigma_{gh}|\nonumber\\
  &\leq & \mathrm{Trace}(\Sigma_{a}^4)^{K_0}\\
  &\leq & C p^{ K_0}. \end{eqnarray*} Then
$\mathbb{E} T_n(Z_1, Z_2)^{K_0}\leq C\frac{p^{K_0}}{n^{4K_0}}.$ Combining this with (\ref{destroy}) and (\ref{dinner}), we get the conclusion.
\end{IEEEproof}

\begin{IEEEproof}[Proof of Lemma \ref{night1}] Let  $g_m= \sum_{1\leq i, j \leq p} \bar{C}^{(\tau)}_ {ij}[( z_{mi}^2-\sigma_{ii})( z_{mj}^2-\sigma_{jj})-2\sigma_{ij}^2]$, then $\mathbb{E} R_1^{K_0}=\frac{1}{n^{2K_0}}(\sum_{m=1}^{n-1}g_m)^{K_0}$.  It is easy to see that  $g_m$ are i.i.d with mean 0,   using moments inequality,  we have $ \mathbb{E}(\sum_{m=1}^{n-1}g_m)^{K_0}\leq L_{K_0}n^{\frac{K_0}{2}}\mathbb{E}(g_m^{K_0}).  $ Now we only need to bound $\mathbb{E}(g_m^{K_0})$. Direct calculation shows that
$
\mathbb{E}(g_m^{K_0})= \sum_{ 1\leq i_1, \cdots, i_{2K_0} \leq p}  \bar{C}^{(\tau)}_{i_1\cdots i_{2K_0}}
\mathbb{E}\prod_{k=1}^{K_0}\{( z_{mi_{2k-1}}^2-\sigma_{i_{2k-1}i_{2k-1}})( z_{mi_{2k}}^2-\sigma_{i_{2k}i_{2k}})-\sigma_{i_{2k-1}i_{2k}}^2\}.
$
Hence,
\begin{eqnarray*}\lbl{hurt}
&&\mathbb{E}(g_m^{K_0})\nonumber \\
&\leq & C \sum_{ 1\leq i_1, \cdots, i_{2K_0} \leq p}\{ \nonumber \\
&&| \mathbb{E}\prod_{k=1}^{K_0} ( z_{mi_{2k-1}}^2-\sigma_{i_{2k-1}i_{2k-1}})( z_{mi_{2k}}^2-\sigma_{i_{2k}i_{2k}})|\nonumber\\
&&+\sum_{l=1}^{K_0}\sigma_{i_{2l-1}i_{2l}}^2 \times \nonumber \\
&&|\mathbb{E}\prod_{k=1, k\neq l}^{K_0}( z_{mi_{2k-1}}^2-\sigma_{i_{2k-1}i_{2k-1}})( z_{mi_{2k}}^2-\sigma_{i_{2k}i_{2k}})|\nonumber\\
&&+\cdots+  \prod_{l=1}^{K_0}C\sigma_{i_{2l-1}i_{2l}}^2\} \nonumber\\
&\leq&  2^{K_0}p^{K_0},
   \end{eqnarray*}
      the last inequality is due to that  there are total $2^{K_0}$ terms  inside the braces and each term can be written as $| \mathbb{E}\prod_{k\in L}(( z_{mi_{2k-1}}^2-\sigma_{i_{2k-1}i_{2k-1}})( z_{mi_{2k}}^2-\sigma_{i_{2k}i_{2k}})|\times \prod_{k \in L^c}\sigma_{i_{2k-1}i_{2k}}^2\}\leq \mathrm{Trace}(\Sigma^2)^{|L_c|}\mathrm{Trace}(\Sigma_a^2)^{|L|}\leq Cp^{K_0}$, where $L\bigcup L_c=\{1, \cdots, K_0\}$. From all above, we have  that
 \begin{eqnarray*}\lbl{hurt}
\mathbb{E}(R_1^{K_0})
&\leq &C\frac{p^{K_0}}{n^{3K_0/2}}.
                         \end{eqnarray*}
            Now define
                  \begin{eqnarray*}  T_k=\sum_{m=2}^{k}H_m,  \,\, T_1=0  \,\, \text{and } \,\,  H_m=\frac{1}{n^2}\sum_{l=1}^{m-1} h(Z_m, Z_l),  \end{eqnarray*}
                  where $h(Z_m, Z_l) = \sum_{1\leq i, j\leq p } 2 b_n\bar{b}_{ij}^{(\tau)}(z_{mi}^2-\sigma_{ii})(z_{lj}^2-\sigma_{jj}).$ Then $R_2=T_{n-1}$.
Let  $\ml{F}_k$ be the $\sigma$-field generated by $(Z_1,\cdots, Z_k)$ for all $1 \leq k$, then $\{T_n, \ml{F}_n, n=1,2,\cdots\}$ is a martingale on the probability space. By the moments bounds for Martingale (Dharmadhikari et al.,1968),  $\mathbb{E}(R_2^{K_0})\leq L_{k_0}n^{\frac{K_0}{2}-1} \sum_{m=2}^{n-1}\mathbb{E}|H_m|^{K_0}$.
Notice $  \mathbb{E} \prod_{j=1}^{K_0}h(Z_m,Z_{l_j})=0$ if exist $l_j\neq l_k\leq m-1$ for any $k\neq j$.  So
\begin{eqnarray*}\lbl{sad}
&&\mathbb{E}(R_2^{K_0})\\
&\leq& \frac{L_{k_0}}{n^{3K_0/2+1}} \sum_{m=2}^{n-1} \sum_{\substack{1\leq l_j\leq m-1\\1\leq j \leq K_0 }}  \mathbb{E} \prod_{j=1}^{K_0}h(Z_m,Z_{l_j})\\
&\leq&  \frac{C}{n^{K_0} }\mathbb{E}h(Z_1,Z_2)^{K_0}.
\end{eqnarray*}  Finally,  we can bound $\mathbb{E}R_2^{K_0}$ by $C\frac{(C_np)^{K_0}}{n^{2K_0}}$, because
           \begin{eqnarray*}\lbl{sad}
&&\mathbb{E}h(Z_1,Z_2)^{K_0}\\
&\leq& (2b_n)^{K_0}C_n^{K_0}\{\sum_{1\leq i_1, \cdots, i_{K_0}\leq p} |\mathbb{E}\prod_{j=1}^{K_0}(z_{1i_j}^2-\sigma_{i_ji_j})|\}^2 \\
&\leq& C\frac{(C_np)^{K_0}}{n^{K_0}}.
          \end{eqnarray*}  So we get the desired moment bound for $R_2$.
          \end{IEEEproof}

\begin{IEEEproof}[Proof of Lemma \ref{fly}] As we argued in the proof of Theorem \ref{Dad}, it is not hard to show $R_c(\tau)$ with banding weight $\omega_{ij}^B(\tau)=I_{(|i-j|<\tau)}$ reaches its minimum at $k_0$.
The tapering weights in Cai et al.(2010) and Cai and Zhou (2010), is defined as $ \omega_{ij}^{C}{(\tau)}= \frac{ \tau - |i-j|}{\lfloor \frac{\tau}{2}\rfloor} $ for  $\lfloor \frac{\tau}{2}\rfloor<|i-j|< \tau$ and $\omega_{ij}^{C}{(\tau)}=\omega_{ij}^B({(\tau)})$ otherwise.  It is easy to check that    $\omega_{ij}^{B}{(\tau)}=\omega_{ij}^{C}{(\tau)}$ for $\tau=1, 2, 3$.  So we only need to consider the bandwidth $k_0\geq 3$.
  Define
     \begin{eqnarray*} &&G_{ij}^{(\tau)}\\
     &=&(\omega_{ij}^{C}(\tau+1) -\omega_{ij}^{C}(\tau))\big\{\omega_{ij}^{C}(\tau)+\omega_{ij}^{C}(\tau+1)-\frac{2n-C_n}{n-1}\big\},\nonumber\\
     &&H_{ij}^{(\tau)}\\
     &=&(\omega_{ij}^{C}(\tau+1) -\omega_{ij}^{C}(\tau))\big\{\omega_{ij}^{C}(\tau)+\omega_{ij}^{C}(\tau+1)+ \frac{n(C_n-2)}{n-1}\big\},    \end{eqnarray*}  while $G_{ij}^{(\tau)}< 0$ and $     \frac{2C_n}{\tau}\geq H_{ij}^{(\tau)}>0$ for all $\lfloor \frac{\tau+1}{2}\rfloor<| i-j|\leq \tau$.
Then   \begin{eqnarray*}
&&R_c(\tau+1)-R_c(\tau)\\
&=&\sum_{\lfloor \frac{
\tau+1}{2}\rfloor<| i-j|\leq  \tau}\frac{n-1}{n} (G_{ij}^{(\tau)}\sigma_{ij}^2+ \frac{1}{n}H_{ij}^{(\tau)} \sigma_{ii}\sigma_{jj}). \end{eqnarray*}
Since $\sigma_{ij}=0$ for any $|i-j|\geq k_0,$ it is easy to see that $R_c(\tau+1)>R_c(\tau)$ for all $\tau\geq 2k_0-3$.
When $\lfloor \frac{\tau+3}{2}\rfloor=| i-j|$, we get that $G_{ij}^{(\tau)}\leq -\frac{1}{C}\frac{1}{\tau^2}$.
Combine this with $k_0\leq o(\sqrt{n/C_n})$, we have\begin{eqnarray*}
&&R_c(\tau+1)-R_c(\tau)\\
&\leq&\sum_{ | i-j|=\lfloor \frac{\tau+3}{2}\rfloor}\frac{n-1}{n} G_{ij}^{(\tau)}\sigma_{ij}^2\\
&&+\sum_{\lfloor \frac{\tau+1}{2}\rfloor<| i-j|\leq  \tau} \frac{n-1}{n^2}H_{ij}^{(\tau)} \sigma_{ii}\sigma_{jj}\\
&\leq&-\frac{b^2p}{C\tau^2}+\frac{CpC_n}{n}<0. \end{eqnarray*}
 for all $ \tau \leq 2k_0-4$.  So we conclude that $R_c(\tau)$ reaches its minimum at $2k_0-3$ for $k_0\geq 3$.
\end{IEEEproof}

\section{Acknowledgment}
THe authors thank the editor, an AE and two referees for their helpful comments.

\begin{IEEEbiographynophoto}{Danning Li}
Danning Li received B.S. (2004) in Mathematics from Northeast Normal University in China and M.S. (2008) in Mathematics from Jilin University in China. She received M.S. and Ph.D. in Statistics from University of Minnesota in 2011 and 2013, respectively. Her research interests include random matrix and high dimensional statistics.
\end{IEEEbiographynophoto}

\begin{IEEEbiographynophoto}{Hui Zou}
Hui Zou received B.S. and M.S. in Physics from University of Science and  Technology of China in 1997 and 1999, respectively.
He received M.S. in Statistics from Iowa State University in 2001 and Ph.D. in Statistics from Stanford University in 2005.
His research interests include high-dimensional inference, statistical learning and computational statistics.
\end{IEEEbiographynophoto}


\begin{thebibliography}{SOSL90}


 \bibitem{Brown}
Brown, B. M. (1971). Martingale central limit theorems. {\it Ann. Math. Statist.} 42, 59--66.

 \bibitem{Bickel1}
Bickel, P. J. and Levina, E. (2008a). Regularized estimation of large covariance matrices. {\it  Ann.  Statist.} 36, 199--227.

 \bibitem{Bickel2}
Bickel, P. J. and Levina, E. (2008b). Covariance regularization by thresholding. {\it Ann. Statist.} 36, 2577--2604.

\bibitem{Bien}
Bien, J., Bunea, F., and Xiao, L. (2015). Convex Banding
of the Covariance Matrix, {\it J. Amer. Statist. Assoc.,} DOI: 10.1080/01621459.2015.1058265

 \bibitem{Cai}
 Cai, T., Zhang, C. H., Zhou, H. (2010). Optimal rates of convergence for covariance matrix estimation. {\it Ann. Statist.} 38, 2118--2144.

\bibitem{CaiLiu}
Cai, T. and Liu, W. (2011). Adaptive thresholding for sparse covariance matrix estimation.
{\it J. Amer. Statist. Assoc.,} 106, 672--684.

\bibitem{CaiYuan}
Cai, T. and Yuan, M. (2012). Adaptive covariance matrix estimation through block thresholding.  {\it Ann. Statist.} 40, 2014--2042.

\bibitem{De}
De la Pe$\tilde{n}$a, V.H. (1999). A general class of exponential inequalities for martingales and ratios. {\it Ann. Probab.}
27, 537--564.

 \bibitem{Dh}
Dharmadhikari, S. W., Fabian, V.,  Jogdeo, K. (1968). Bounds on the moments of martingales. {\it The Annals of Mathematical Statistics} 1719--1723.

 \bibitem{DJ}
Donoho, D. and Johnstone, I. (1995). Adapting to unknown smoothness via wavelet shrinkage.  {\it J. Amer. Statist. Assoc.} 90, 1200--1224.

  \bibitem{E1}
  Efron, B. (1986). How biased is the apparent error rate of a prediction rule. {\it J. Amer. Statist. Assoc.,} 81, 461--470.

    \bibitem{E2}
   Efron, B. (2004). The estimation of prediction error: covariance penalties and cross-validation. {\it J. Amer. Statist. Assoc.} 99, 619--632.

\bibitem{E3}
Efron, B., Hastie, T., Johnstone, I. and Tibshirani, R. (2004). Least angle regression (with discussion). {\it Ann. Statist.}  32, 407--499.

\bibitem{EKaroui}
El Karoui, N. (2008). Operator norm consistent estimation of large dimensional sparse co- variance matrices. {\it Ann.  Statist.} 36, 2717--2756.

\bibitem{Furrer}
   Furrer, R. and Bengtsson, T. (2007). Estimation of high-dimensional prior and posterior covariance matrices in Kalman filter variants. {\it J. Multivariate Anal.}  98(2), 227--255.

\bibitem{Hall}
 Hall, P. (1984).  Central limit theorem for integrated square error of multivariate nonparametric density estimators. {\it J. Multivariate Anal.} 14, 1--16.

 \bibitem{John}
Johnstone, I. M. (2001). On the distribution of the largest eigenvalue in principal components analysis. {\it Ann.  Statist.} 29, 295--327.

  \bibitem{Murri}
Muirhead, R. J. (1982). Aspects of Multivariate Statistical Theory.  Wiley, New York.

  \bibitem{Qiu}
Qiu, Y. and Chen, S. (2012). Test for bandedness of high-dimensional covariance matrices and bandwidth
estimation. {\it Ann. Statist.}, 40, 1285--1314.


   \bibitem{Rothman}
Rothman, A. J., Levina, E. and Zhu, J. (2009). Generalized thresholding of large covariance matrices. {\it J. Amer. Statist. Assoc.} 104(485), 177--186.

   \bibitem{Shao97}
 Shao, J. (1997).  An asymptotic theory for linear model selection (with discussion). {\it Statistica Sinica}, 7, 221--242.

 \bibitem{stein}
Stein, C. (1981). Estimation of the mean of a multivariate normal distribution. {\it Ann. Statist.} 9(6), 1135--1151.

  \bibitem{Withers}
Withers, C.S. (1985) The moments of the multivariate normal, {\it Bulletin of the Australian Mathematical Society} 32, 103--107.

 \bibitem{Wu}
 Wu, W. B. and Pourahmadi, M. (2003). Nonparametric estimation of large covariance matrices of longitudinal data. {\it Biometrika} 90(4), 831--844.

 \bibitem{Xiao}
Xiao, L. and Bunea, F. (2014). On the theoretic and practical merits of the banding estimator for
large covariance matrices. {\it arXiv preprint, arXiv:1402.0844.}

 \bibitem{yang}
Yang, Y. (2005). Can the strengths of AIC and BIC be shared? A conflict between model indentification and regression estimation. {\it Biometrika} 92(4), 937--950.

\bibitem{Feng}
 Yi, F. and Zou, H. (2013). SURE-tuned tapering estimation of large covariance matrices. {\it Computational Statisticsa and Data Analysis,} 58, 339--351.

\bibitem{Zou07}
Zou, H., Hastie, T. and Tibshirani, R. (2007). On the degrees of freedom of the lasso. {\it Ann. Statist.} 35, 2173--2192.


 \end{thebibliography}
\end{document}